\documentclass[11pt,reqno]{amsart}

\usepackage[T1]{fontenc}
\usepackage[english]{babel}
\usepackage[a4paper,margin=2.7cm]{geometry}
\usepackage{microtype}
\usepackage{dsfont}
\usepackage{amssymb}
\usepackage{mathtools}
\usepackage{amsthm}
\usepackage{graphicx}
\usepackage{xcolor}
\usepackage[square,numbers,comma,sort&compress]{natbib}
\usepackage{hyperref}
\hypersetup{colorlinks=true,linkcolor=blue,urlcolor=blue,citecolor=blue}
\definecolor{fqgreen}{rgb}{0,0.45,0}
\allowdisplaybreaks

\usepackage[normalem]{ulem}

\newcommand{\Lop}{\mathcal L}
\newtheorem{theorem}{Theorem}[section]
\newtheorem{proposition}[theorem]{Proposition}
\newtheorem{lemma}[theorem]{Lemma}
\newtheorem{cor}[theorem]{Corollary}
\theoremstyle{definition}
\newtheorem{definition}[theorem]{Definition}
\theoremstyle{remark}
\newtheorem{remark}[theorem]{Remark}
\numberwithin{equation}{section}

\DeclareMathOperator*{\esssup}{ess\,sup}
\DeclareMathOperator*{\esslim}{ess\,lim}
\DeclareMathOperator*{\esslimsup}{ess\,lim\,sup}
\newcommand{\sgn}{\operatorname{sgn}}
\newcommand{\eps}{\varepsilon}
\newcommand{\rr}{\mathbb{R}}
\newcommand{\nn}{\mathbb{N}}

\title[Convection--diffusion equations with non-integrable data]{Convection--diffusion equations with non-integrable initial data}

\author[N.~De Nitti]{Nicola De Nitti}
\address[N.~De Nitti]{Politecnico di Bari, Dipartimento di Meccanica, Matematica e Management, Via E.~Orabona 4, 70125 Bari, Italy.}
\email{nicola.denitti\@@{}poliba.it}

\author[L.~I.~Ignat]{Liviu I. Ignat}
\address[L.~I.~Ignat]{Institute of Mathematics ``Simion Stoilow'' of the Romanian Academy, 21 Calea Grivitei Street, 010702 Bucharest, Romania.
\newline\indent
Academy of Romanian Scientists, Ilfov Street, no.~3, Bucharest, Romania.}
\email{liviu.ignat\@@{}gmail.com}

\author[F.~Quir\'{o}s]{Fernando Quir\'{o}s}
\address[F.~Quir\'{o}s]{Departamento de Matem\'aticas, Universidad Aut\'onoma de Madrid (UAM), Campus de Cantoblanco, 28049 Madrid, Spain,
\newline\indent
and Instituto de Ciencias Matem\'aticas
\newline\indent
ICMAT (CSIC-UAM-UCM-UC3M), 28049 Madrid, Spain.}
\email{fernando.quiros\@@{}uam.es}

\begin{document}

\begin{abstract}
We characterize the long-time behaviour of solutions to the convection--diffusion equation $u_t-u_{xx}+(|u|^{q-1}u)_x=0$, $q>1$, with initial data behaving like a non-integrable power at infinity, $|x|^\gamma u_0(x)\to k_\pm$ as $x\to\pm\infty$, $0<\gamma\le1$. The decay rate in time and the asymptotic profile depend on the strength of convection as compared to the size of the solution at infinity, leading either to a convective behaviour or to a diffusive one. When the nonlinear convection dominates, we observe an enhanced dissipation effect: the decay rate in time is faster than for the purely diffusive equation with the same initial datum. For $0<\gamma<1$ we obtain the limit self-similar profile, which, in the convective regime, has one shock, or two when $k_+k_-<0$ and $|k_+|>\eps(\gamma,q)|k_-|$, where $\eps(\gamma,q)\in(0,1)$ is an explicit threshold. For $\gamma=1$  logarithmic corrections appear, and the limit is a multiple of the heat kernel if $q>2$ and an $N$-wave if $1<q\le2$. The convergence takes place in the spaces to which the datum belongs, and not only on expanding sets. The main new ingredients are an Ole\u{\i}nik type estimate of H\"older type valid for $q\geq 2$ and  a far-field estimate, uniform in time,  showing that beyond the self-similar region the solution is still described by its own datum.
\end{abstract}

\keywords{Convection--diffusion equations, non-integrable initial data, large-time behaviour, self-similar solutions, entropy solutions, uniqueness, one-sided estimates, Marcinkiewicz spaces}
\subjclass[2020]{35B40, 35K55, 35L65, 35L67, 35A02}

\maketitle

\section{Introduction and main results}\label{sect-Introduction}

We study the large-time behaviour of solutions to the Cauchy problem
\begin{equation}
\label{eq.main}
	u_t-u_{xx}+(F(u))_x =0 \quad\textup{in }Q:=(0,\infty)\times \rr,\qquad u(0)=u_0 \quad\textup{on }\rr,
\end{equation}
where
\begin{equation*}
    F(u):= a \, {|u|^{q-1}u},\quad  a\in \{\pm 1\}, \quad q>1,
\end{equation*}
  for bounded, possibly sign-changing, initial data that decay at infinity like a non-integrable power, namely
\begin{equation}\label{initial.gamma}
u_0\in L^\infty(\rr),\qquad \esslim_{x\to\pm\infty}|x|^\gamma u_0(x)=k_\pm,\qquad 0<\gamma\le1, \ (k_+,k_-)\neq (0,0).
\end{equation}
For clarity of the presentation we consider the case $a=1$.

\subsection{Precedents}
When $u_0\in L^1(\rr)$ has mass $M=\int_\rr u_0$, the large-time behaviour of \eqref{eq.main} is well understood and depends on $q$. For $q>2$ convection is negligible and the solution behaves like $MG_t$, $G_t$ being the heat kernel, and for $q=2$ it approaches the source-type solution of \eqref{eq.main} with mass $M$ \cite{EZ}; for $1<q<2$ the hyperbolic part dominates, and suitably rescaled solutions converge to the $N$-wave, that is, to the entropy solution of $u_t+(F(u))_x=0$ with initial datum $M\delta_0$ \cite{EVZArma,Carpio1996,MR735207}. We refer to \cite{EVZIndiana,EsZu97} for several space dimensions, to \cite{LaurencotFast} for $0<q<1$, and to \cite{MR3789847,endalJMPA,endalDCDS} for nonlocal diffusions. For bounded non-integrable data the large-time behaviour of evolution equations may be considerably more involved \cite{zbMATH01799645}. For the heat equation and the porous medium equation in $\rr^N$, without convection, data decaying like $|x|^{-\gamma}$ with $0<\gamma<N$ lead to self-similar asymptotics, while the borderline case $\gamma=N$ produces logarithmic corrections; see \cite{AlRo84,KaUg87,Herraiz-1999}.

The case $q=2$, $0<\gamma<1$, with nonnegative data, for which \eqref{eq.main} is the viscous Burgers equation up to a constant factor in $u$, was considered by Ghoul, Masmoudi and Pacherie \cite{MR4687278,MR4899801}. They  proved that nonlinear convection produces an \emph{enhanced dissipation} effect: the solution decays like $t^{-\gamma/(1+\gamma)}$, faster than the solution $G_t\ast u_0$ of the heat equation, which decays like $t^{-\gamma/2}$. They also described the asymptotic profile, which has one discontinuity, for $k_-=\pm k_+\neq 0$ \cite[Theorem~1.5, Proposition~1.7]{MR4687278}; see also \cite{2406.10874} for a refined description near the discontinuity.  For $\gamma=1$, data with $k_+=-k_-$ were considered in \cite{Dix2002}.

Related results for Burgers-type equations with fractional or critical dissipation and singular or non-decaying data can be found in \cite{JaSe23,AlBe22,KarchMiaoXu}.

\subsection{Heuristics and scope of the results}
The decay rate $t^{-\gamma/(1+\gamma)}$ found in \cite{MR4687278} for $q=2$ is one instance of a general competition between diffusion and convection, which decides the large-time behaviour of \eqref{eq.main}--\eqref{initial.gamma} and is best understood by looking first at the two extreme cases. Without convection, the solution $G_t\ast u_0$ of the heat equation with a datum decaying like $|x|^{-\gamma}$, $\gamma<1$, decays like $t^{-\gamma/2}$ and lives on the diffusive scale $|x|\sim t^{1/2}$. Without diffusion, the entropy solution of the conservation law
\begin{equation}\label{conservation.law}
    u_t+(F(u))_x=0\quad\text{in }Q,\qquad u(0)=u_0\quad\text{on }\rr,
\end{equation}
with the same datum decays like $t^{-\alpha}$ and lives on the scale $|x|\sim t^\beta$, where
\begin{equation}\label{alpha.beta}
    \alpha:=\frac{\gamma}{1+(q-1)\gamma},\qquad\beta:=\frac1{1+(q-1)\gamma} .
\end{equation}
Indeed, the hyperbolic scaling $u_\lambda(t,x):=\lambda^\gamma u(\lambda^{1/\beta}t,\lambda x)$ preserves both the conservation law and the behaviour of the datum at infinity, and the convergence of $u_\lambda$ as $\lambda\to\infty$ translates into large-time behaviour, since for $\lambda=t^\beta$, $u_\lambda(1,\xi)=t^{\alpha}u(t,t^\beta\xi)$. This suggests that $u(t,x)\approx t^{-\alpha}V(xt^{-\beta})$ as $t\to\infty$, with $V$ the limit of $u_\lambda(1,\cdot)$. We shall prove that the solutions of \eqref{eq.main} satisfy both decay estimates, the parabolic one and the hyperbolic one. Hence they decay at the faster of the two rates, $t^{-\max\{\gamma/2,\,\alpha\}}$. Since $\alpha>\gamma/2$ if and only if $q<1+ 1/\gamma$, this gives the critical curve $q=1+1/\gamma$. If $q$ is below this critical value,  the solution decays faster than the solution of the heat equation with the same datum: this is the enhanced dissipation of \cite{MR4687278}, whose rate $t^{-\gamma/(1+\gamma)}$ is $t^{-\alpha}$ for $q=2$.

The same competition determines the asymptotic profile. Above the curve, $q>1+1/\gamma$, the decay is the parabolic one, and Duhamel’s formula shows that convection is a lower-order perturbation, since $u(t)-G_t\ast u_0$ is $o(t^{-\gamma/2})$; hence the solution converges to that of the heat equation with datum $U_0(x):=k_{\sgn x}|x|^{-\gamma}$, the limit of the data under the scaling $\lambda^\gamma u_0(\lambda x)$. Below the curve, $q<1+\frac1\gamma$, the decay is the hyperbolic one, and we use the hyperbolic scaling $u_\lambda$, which transforms \eqref{eq.main} into
\[
    \partial_tu_\lambda+\partial_x\big(F(u_\lambda)\big)=\lambda^{(q-1)\gamma-1}\,\partial_{xx}u_\lambda ,
\]
a conservation law with vanishing viscosity. Formally, and we shall make this rigorous, $u_\lambda$ converges to the entropy solution of \eqref{conservation.law} with datum $U_0$. On the curve $q=1+\frac1\gamma$ the two rates coincide, the viscosity in the rescaled equation is equal to one, and the limit is the solution of \eqref{eq.main} itself with datum $U_0$. In the three cases the limit is self-similar.

The borderline case $\gamma=1$ is different. The scaling that preserves the datum is now $\lambda u_0(\lambda x)$, and its limit $U_0=k_\pm|x|^{-1}$ is not locally integrable at the origin. This reflects the fact that the mass of $u_0$ in large balls diverges logarithmically, 
\begin{equation}
  \label{eq:mass.R}
\int_{|x|<R}u_0\,dx=(k_++k_-)\log R+o(\log R),  
\end{equation}
 and that the scaling brings these balls to the origin. The decay rates acquire a logarithmic factor.    When $k_++k_-\neq 0$, the asymptotic profiles are, as for integrable data with $q\neq 2$, the heat kernel if $q>2$ and the $N$-wave if $q<2$,  but with a mass that grows like $\log t$. However, when $q=2$, the profile here is not given by the fundamental solutions of   the complete equation, as in the integrable case,  but by the $N$-wave with corrected mass. 

\subsection{Main results}
Throughout the paper $u$ denotes the unique solution of \eqref{eq.main} given by Theorem~\ref{integrable} below, $G_t$ is the heat kernel, $U_0(x)=k_{\sgn x}|x|^{-\gamma}$, $\alpha$ and $\beta$ are as in \eqref{alpha.beta}, and $\theta:=\max\{\frac12,\beta\}$. Thus $\theta=\beta>\frac12$ if and only if $q<1+\frac1\gamma$. The scale $t^\theta$ is the larger of the diffusive scale $t^{1/2}$ and the convective one $t^\beta$, and $\gamma\theta=\max\{\frac\gamma2,\alpha\}$.

We treat the two cases $0<\gamma<1$ and $\gamma=1$ separately.

\begin{theorem}\label{first.term.gamma<1}
Let $\gamma\in(0,1)$, $q>1$, $k_\pm\in\rr$, and  $u_0$ satisfying \eqref{initial.gamma}. Then there is a constant $C$ such that
\begin{equation}\label{thm1.decay}
    \|u(t)\|_{L^\infty(\rr)}\le C(1+t)^{-\gamma\theta}\qquad\text{for all }t\ge0.
\end{equation}
Moreover,
\begin{equation}\label{thm1.conv}
    \lim_{t\to\infty}\big\|u(t)-U(t)\big\|_{L^{1/\gamma,\infty}(\rr)}=0 ,
\end{equation}
and, for every $p\in(\frac1\gamma,\infty)$, ($p\in (\frac 1\gamma,\infty]$ if $q\geq 1+\frac 1\gamma$),
\begin{equation}\label{thm1.conv.p}
    \lim_{t\to\infty}t^{\theta(\gamma-\frac1p)}\big\|u(t)-U(t)\big\|_{L^p(\rr)}=0 ,
\end{equation}
where $U$ is
\begin{enumerate}
\item[(a)] the solution $U(t)=G_t\ast U_0$ of the heat equation with initial datum $U_0$, if $q>1+\frac1\gamma$;
\item[(b)] the unique solution of \eqref{eq.main} with initial datum $U_0$ such that
    \begin{equation}\label{eq:uniqueness.class}
        \sup_{t>0}t^{\gamma\theta}\|U(t)\|_{L^\infty(\rr)}<\infty,
    \end{equation}
    if $q=1+\frac1\gamma$;
\item[(c)] the unique entropy solution of \eqref{conservation.law} with initial datum $U_0$ such that \eqref{eq:uniqueness.class} holds, if $1<q<1+\frac1\gamma$; for $(k_+,k_-)\neq(0,0)$ it is described in Theorem~\ref{self.profiles}.
\end{enumerate}
In all three cases $U$ is self-similar, $U(t,x)=t^{-\gamma\theta}U_\gamma(xt^{-\theta})$ with $U_\gamma:=U(1,\cdot)$.
\end{theorem}

Data satisfying \eqref{initial.gamma} belong to $L^{1/\gamma,\infty}(\rr)\cap L^\infty(\rr)$ and to $L^p(\rr)$ for every $p>\frac1\gamma$, and, if $(k_+,k_-)\neq(0,0)$, to no smaller Lebesgue or Marcinkiewicz space; the convergence \eqref{thm1.conv}--\eqref{thm1.conv.p} therefore takes place in the spaces to which the datum belongs. When $(k_+,k_-)=(0,0)$ the result still holds but the decay \eqref{thm1.decay} is not optimal and $\|u(t)\|_{L^{1/\gamma,\infty}(\rr)}\rightarrow 0$ as $t\rightarrow \infty$ without a profile. 

Since $U_0$ is unbounded, the problems in (b) and (c) are outside the classical well-posedness theories, and solutions are understood in the sense of Definition~\ref{def.entropy}. They take the datum in $L^1_{\textup{loc}}(\rr)$. To prove uniqueness we   require some control on the solution as $t\downarrow0$.  The condition $\sup_{t>0}t^{\gamma\theta}\|U(t)\|_{L^\infty(\rr)}<\infty$ is  natural: it is the decay of the self-similar solutions, and, since $\gamma\theta<\frac1q$ when $q\leq 1+1/\gamma$, it falls within the uniqueness result of Theorem~\ref{unique.entropy}.

If $(k_+,k_-)\neq(0,0)$, estimates \eqref{thm1.conv}--\eqref{thm1.conv.p} say that, for $X=L^{1/\gamma,\infty}(\rr)$, for $X=L^p(\rr)$, $\frac1\gamma<p<\infty$, and, if $q\ge1+\frac1\gamma$, for $X=L^\infty(\rr)$,
\[
    \|u(t)-U(t)\|_X=o\big(\|U(t)\|_X\big)
\]
or equivalently,
\[
    \|t^{\gamma\theta}u(t,t^\theta\cdot)-U_\gamma\|_X\to0.
\]

In the convective regime $1<q<1+\frac1\gamma$ the endpoint $p=\infty$ cannot be included when $(k_+,k_-)\neq(0,0)$, since $u(t)$ is continuous while $U_\gamma$ has jumps.

Convergence in $L^1_{\textup{loc}}(\rr)$ of the rescaled solutions, the form in which such results are usually stated, is much weaker, since it says nothing outside the region $|x|\lesssim t^\theta$. What happens there is the content of our next result, which is also the key to \eqref{thm1.conv}--\eqref{thm1.conv.p}.

\begin{theorem}\label{farfield}
    Let $\gamma\in(0,1)$, $q>1$, $k_\pm\in\rr$, and  $u_0$ satisfying \eqref{initial.gamma}. Then
    \begin{equation}\label{thm.ff}
        \lim_{L\to\infty}\ \sup_{t>0}\ \sup_{\pm x\ge L(1+t^\theta)}\big|\,|x|^\gamma u(t,x)-k_\pm\big|=0 ,
    \end{equation}
    and, for some $\mathcal{N}<\infty$,
    \begin{equation}\label{thm.env}
        \big(1+t^\theta+|x|\big)^\gamma|u(t,x)|\leq \mathcal{N}\quad \text{for all}\ t>0,\  x\in \rr.
    \end{equation}
\end{theorem}

In the far field $|x|\gg t^\theta$ the solution is thus still described by its own datum, uniformly in time, and not by the self-similar profile; the two descriptions match, since $|\xi|^\gamma U_\gamma(\xi)\to k_\pm$ as $\xi\to\pm\infty$.

\begin{theorem}\label{first.term.gamma=1}
    Let $\gamma=1$, $q>1$,  $k_\pm\in\rr$ and  $u_0$ satisfying \eqref{initial.gamma}. Then there is a constant $C$ such that
    \begin{equation}\label{thm2.decay}
        \|u(t)\|_{L^\infty(\rr)}\le
        \begin{cases}
         C\,\frac{\log(2+t)}{(1+t)^{1/2}}&\quad\text{if }q>2,\\
     C\Big(\frac{\log(2+t)}{1+t}\Big)^{1/q}&\quad\text{if }1<q\le2      .
        \end{cases}
    \end{equation}
    Moreover,
    \begin{equation}\label{thm2.weak}
        \lim_{t\to\infty}\frac1{\log t}\,\big\|u(t)-\mathcal U(t)\big\|_{L^{1,\infty}(\rr)}=0,
    \end{equation}
    and
         \begin{equation}
        \label{thm2.heat}
        \lim_{t\to\infty}\frac{t^{\frac12(1-\frac1p)}}{\log t}\,\big\|u(t)-\mathcal U(t)\big\|_{L^p(\rr)}=0\quad \text{if } q>2,\ p\in (1,\infty],
        \end{equation}
        \begin{equation}
        \label{thm2.nwave}
        \lim_{t\to\infty}\frac{t^{\frac{p-1}{pq}}}{(\log t)^{\frac{p+q-1}{pq}}}\,\big\|u(t)-\mathcal U(t)\big\|_{L^p(\rr)}=0\quad \text{if } 1<q\le2,\ p\in (1,\infty),
    \end{equation}
    where
    \[
        \mathcal U(t):=\frac{k_++k_-}2\,\log t\;G_t\quad\text{if }q>2,\qquad \mathcal U(t):=U_{\frac{k_++k_-}q\log t}(t)\quad\text{if }1<q\le2,
    \]
    and $U_M$ is the $N$-wave  of mass   $M\neq0$,
    \begin{equation}\label{nwave}
        U_M(t,x)=\sgn(M)\Big(\frac{x}{qt}\Big)^{\frac1{q-1}}\mathds{1} _{\{0<x<q(|M|/(q-1))^{(q-1)/q}\,t^{1/q}\}}.
    \end{equation}
    When $M=0$, $\mathcal{U}\equiv 0$.
\end{theorem}

\begin{remark}The above $U_M$ is the $N$-wave of mass $M$, that is, the entropy solution of \eqref{conservation.law} with initial datum $M\delta_0$, which is unique among the entropy solutions in $L^\infty((0,\infty);L^1(\rr))\cap L^\infty((\tau,\infty)\times\rr)$, $\tau>0$ \cite{MR735207}. As in the case $\gamma\in(0,1)$, the endpoint $p=\infty$ is excluded for $1<q\le2$, since $u(t)$ is continuous while the $N$-wave has a jump.

The limit profile, $\mathcal U(t)$, is the heat kernel of mass $\frac{k_++k_-}2\log t$ if $q>2$, and the $N$-wave of mass $\frac{k_++k_-}q\log t$ if $1<q\le2$.
 In contrast with the case of integrable data, its mass grows with time.
It  is, up to lower-order terms, the mass of the datum in the region $|x|\lesssim t^\theta$ which grows logarithmically, see \eqref{eq:mass.R}.

  When $k_++k_-=0$ the decay is faster  than in \eqref{thm2.decay}. To obtain the precise decay one needs to have more information  on the behaviour of the initial data at infinity, see \cite{Dix2002} for Burgers' equation.
\end{remark}

The far-field description of Theorem~\ref{farfield} also holds for $\gamma=1$, at the price of a logarithm in the scale at which the far field begins. Let
\begin{equation}\label{m.def}
    m(t):=\begin{cases}(1+t)^{1/2},&q>2,\\[2pt] (1+t)^{1/q}\big(\log(2+t)\big)^{\frac{q-1}q},&1<q\le2,\end{cases}
\end{equation}
which is, up to constants, the length of the region where $\mathcal U(t)$ lives.

\begin{theorem}\label{farfield.gamma1}
    Let $\gamma=1$,  $q>1$,  $k_\pm\in\rr$, and  $u_0$
 satisfying
 \eqref{initial.gamma}. Then
    \begin{equation}\label{ff.gamma1}
        \lim_{L\to\infty}\ \sup_{t>0}\ \sup_{\pm x\ge L\,m(t)\log(2+t)}\big|\,|x|\,u(t,x)-k_\pm\big|=0 ,
    \end{equation}
    and, for some $\mathcal{M}<\infty$,
    \begin{equation}\label{env.gamma1}
   \frac{x^2\,|u(t,x)|}{|x|+m(t)\log(2+t)}\leq \mathcal{M}\quad \text{for all}\ t>0,\ x\in \rr.
    \end{equation}
\end{theorem}

The global bound \eqref{env.gamma1} connects the decay $\log(t)/m(t)$ in inner regions, $|x|\lesssim m(t)$,  with the decay, $1/(m(t)\log (t))$, in the far field, $|x|\gtrsim m(t)\log(t)$, given by \eqref{ff.gamma1}. It is all we need in the proof of Theorem \ref{first.term.gamma=1}.

The profiles in Theorem~\ref{first.term.gamma<1}\,(c) are described by our next result. Let $r(q)$ be the unique root in $(0,1)$ of
\begin{equation}\label{r.q}
    (q-1)r^q+qr^{q-1}=1
\end{equation}
and set
\begin{equation*}
    \eps(\gamma,q):=r(q)\Big(\frac{\gamma}{1+r(q)-\gamma}\Big)^{\gamma}\in(0,1).
\end{equation*}

\begin{theorem}\label{self.profiles}
    Let  $\gamma\in(0,1)$, $q>1$,  $k_\pm\in\rr$. There exists a unique entropy solution $U$ of \eqref{conservation.law} with initial datum $U_0$, in the sense of Definition~\ref{def.entropy}, such that
    \[
        \sup_{t>0}t^\alpha\|U(t)\|_{L^\infty(\rr)}<\infty.
    \]
    It is self-similar, $U(t,x)=t^{-\alpha}V(xt^{-\beta})$, where $V\in L^\infty(\rr)$ is smooth outside a finite set, at each point of which it has one-sided limits, and $|\xi|^\gamma V(\xi)\to k_\pm$ as $\xi\to\pm\infty$. Moreover:
    \begin{enumerate}
        \item[(i)] if $k_+k_->0$, then $V$ has exactly one discontinuity;
        \item[(ii)] if $k_+k_-<0$ and $|k_+|\le\eps(\gamma,q)|k_-|$, then $V$ has exactly one discontinuity;
        \item[(iii)] if $k_+k_-<0$ and $|k_+|>\eps(\gamma,q)|k_-|$, then $V$ has exactly two discontinuities, and $|V|$ is increasing between them;
        \item[(iv)] if exactly one of $k_\pm$ vanishes, then $V$ has exactly one discontinuity; moreover $V$ vanishes to the right of it if $k_+=0$, and $V$ vanishes on $(-\infty,0]$ and $|V(\xi)|=\big(\frac\xi q\big)^{\frac1{q-1}}$ between the origin and the discontinuity if $k_-=0$.
    \end{enumerate}
\end{theorem}

The number $r(q)$ is the tangency ratio for the flux $F$: the chord of $F$ joining $-1$ and $\rho>0$ is tangent to the graph of $F$ at $\rho$ if and only if $\rho=r(q)$. In case (iii) the first shock joins a left state $v_l$ and a right state $v_r$ of opposite signs with $|v_r|/|v_l|=r(q)$, so that, by homogeneity, its speed, the slope of the chord of $F$ joining $v_l$ and $v_r$, equals the characteristic speed $F'(v_r)$ of its right state. Since $\eps(\gamma,q)<1$, the symmetric case $k_+=-k_-$ always leads to two shocks. In case (iv) with $k_-=0$ the profile contains the branch $V(\xi)=\sgn(k_+)\big(\frac\xi q\big)^{\frac1{q-1}}$, which corresponds to $U(t,x)=\sgn(k_+)\big(\frac x{qt}\big)^{\frac1{q-1}}$, the branch that also appears in the $N$-wave \eqref{nwave}; the profile is continuous, but not smooth, at the origin. See Figure~\ref{fig:profiles} in Section \ref{sect-profiles}.

For the Burgers flux and data with $k_-<0<k_+$ the structure of the profile was left open in \cite{MR4687278}. We stress that, for sign-changing solutions, the odd flux $F(u)=|u|u$ considered here is \emph{not} the Burgers flux, since it is convex on $[0,\infty)$ and concave on $(-\infty,0]$. This change of convexity is responsible for the two-shock profiles of Theorem~\ref{self.profiles} above.
Even fluxes, such as $|u|^q$, lead to a different behaviour, since for integrable sign-changing data the limits are $N$-waves determined by two parameters \cite{MR735207,Kim2004}. This will not be treated here.

\subsection{Strategy}
The proofs of Theorems~\ref{first.term.gamma<1} and \ref{first.term.gamma=1} start with decay estimates, which hold in all cases. After them, the argument depends on the regime. If diffusion dominates, that is, if $q>1+\frac1\gamma$, Duhamel's formula and the decay estimates show directly that the solution is close to the solution of the heat equation. In the other cases we follow the scaling approach of \cite{EZ,EVZArma}: compactness of the rescaled solutions, and identification of the limit through a uniqueness result. Each step requires new ingredients.

\noindent\emph{Decay.} In Section~\ref{sect-estimates} we show, for every flux $F\in C^1(\rr)$ with $F(0)=0$, that $\int_0^su(t)^*\le\int_0^su_0^*$ for all $s,t>0$, where $f^*$ is the decreasing rearrangement of $|f|$: the integral of $|u(t)|$ over a set of measure $s$ never exceeds the largest integral of $|u_0|$ over a set of measure $s$. This gives bounds in $L^{p,\infty}(\rr)$, $p>1$, an $L^{p_0,\infty}$--$L^\infty$ decay estimate (Theorem~\ref{week.integrable}) and, combined with one-sided estimates, the optimal decay rates in the convective regime (Theorem~\ref{th.decay.p_0}).
For the case $p_0=1$ we need some control of the tails as in \eqref{initial.gamma}.
 For $q\ge2$ we prove, by a two-point maximum principle, the one-sided H\"older estimate
\[u(t,x)-u(t,y)\le\Big(\frac{x-y}{qt}\Big)^{1/(q-1)},\quad  y<x,\] for nonnegative solutions (Theorem~\ref{holder-oleinik}); it is uniform in the viscosity, and the $N$-waves \eqref{nwave} with $M>0$ attain equality in it.

\noindent\emph{The diffusive regime.} If $\gamma<1$ and $q>1+\frac1\gamma$, the decay estimates and Duhamel's formula give $u(t)-G_t\ast u_0=o(t^{-\gamma/2})$ in $L^\infty(\rr)$ (Subsection~\ref{subsect.supercritical}), and the limit is the solution of the heat equation with datum $U_0$. For $\gamma=1$ and $q>2$ the same argument gives
\[\frac{t^{1/2}}{\log t}\|u(t)-G_t\ast u_0\|_{L^\infty(\rr)}\to 0\](Subsection~\ref{subsect.heat1}), and the limit is the heat kernel, with a mass that grows like $\log t$.

\noindent\emph{Compactness.}  If $\gamma<1$ and $q\le 1+\frac1\gamma$, we rescale. By the decay estimates, the rescaled solutions are bounded in $L^\infty((\tau,\infty)\times\rr)$ for every $\tau>0$, and they solve a viscous conservation law with viscosity $\eps_\lambda=\lambda^{(q-1)\gamma-1}\le1$, and a local energy estimate shows that $\eps_\lambda^{1/2}\partial_xu_\lambda$ is bounded in $L^2_{\textup{loc}}$. When $q=1+\frac1\gamma$, where $\eps_\lambda=1$, this gives compactness in $L^1_{\textup{loc}}$ by the Aubin--Lions--Simon lemma. When $q<1+\frac1\gamma$ the viscosity vanishes in the limit, and compactness follows from the compensated compactness theory of Tartar and Murat, which only requires $F$ not to be affine on any interval.

\noindent\emph{Uniqueness for the limit problems.} The limit problems have unbounded data, and the classical uniqueness theory does not apply. In Section~\ref{sect-CL} we prove uniqueness, for the conservation law and for \eqref{eq.main}, for every locally integrable datum, among the solutions with $\|u(t)\|_{L^\infty(\rr)}\lesssim t^{-\alpha}$ as $t\downarrow0$, for any $\alpha<\frac1q$ (Theorem~\ref{unique.entropy}), by a simple argument, based on a recursive inequality between balls of radii $R$ and $2R$ and an iteration.
For $1<q\le2$, the limit problem of Theorem~\ref{first.term.gamma=1} has a Dirac mass as datum, and we rely on \cite{MR735207}.

\noindent\emph{Identification of the limit.} That the limit $U$ takes the datum $U_0$ is proved by a tail estimate, based on the time derivative of $\int_\rr|u_\lambda|\phi\,dx$ for compactly supported $\phi$. Together with the decay $\|U(t)\|_{L^\infty(\rr)}\lesssim t^{-\alpha}$, $\alpha<\frac1q$, this places $U$ in the uniqueness class of Theorem~\ref{unique.entropy}. Uniqueness shows that the limit does not depend on the sequence, so that the whole family converges. Since $U_0$ and the uniqueness class are invariant under the scaling,  $U$ is self-similar.

\noindent\emph{Convergence at every fixed time.} Compactness gives convergence of the rescaled solutions in $L^1_{\textup{loc}}$ of $Q=(0,\infty)\times\rr$, that is, for almost every time, while the theorems are statements at any fixed time. To pass from almost every time to any time we prove a fixed-time compactness principle, based on a weighted $L^1$ stability estimate for the viscous equation (Appendix~\ref{sect-fixedtime}).

\noindent\emph{The far field.} The convergence \eqref{thm1.conv}--\eqref{thm1.conv.p} and \eqref{thm2.weak}--\eqref{thm2.nwave}, in norms that see the whole line, requires a description of the solution beyond the self-similar region, which the rescalings do not reach. We obtain it in Theorems~\ref{farfield} and \ref{farfield.gamma1} by comparison on moving domains of the form $\{|x|>\Lambda m(t)\}$, with an explicit barrier containing a term that controls the solution on the moving boundary; beyond the decay estimates \eqref{thm1.decay} and \eqref{thm2.decay}, no gradient bound, one-sided estimate, convexity of $F$ or sign condition is used, so that all $q>1$ are covered at once. Descriptions in which an inner self-similar region is matched with an outer region where the datum persists appear for parabolic problems with slowly decaying data \cite{Herraiz-1999}; for \eqref{eq.main} we are not aware of a result of this type. 

\noindent\emph{The case $\gamma=1$.} The mass of $\lambda u_0(\lambda x)$ in $\{|x|<1\}$ is $\approx(k_++k_-)\log\lambda$, so the rescalings carry a logarithmic correction, with which the rescaled data converge to the Dirac mass $(k_++k_-)\delta_0$ and the limit problem is that of integrable data. For $q>2$, as in the diffusive regime for $\gamma<1$, Duhamel's formula and the decay estimates give the heat kernel. For $1<q\le2$ we show that the limit $U$ of the solutions rescaled with the logarithmic correction is integrable, with mass $k_++k_-$, and takes the initial datum $(k_++k_-)\delta_0$, where the tails are controlled by comparison with the solution with datum $|u_0|$, which is nonnegative, and the mass by the weak formulation with cut-off functions. The uniqueness theorem of Liu and Pierre \cite{MR735207} then identifies $U$ as the $N$-wave.

\noindent\emph{Self-similar profiles.} The profile $V$ of a self-similar solution $U(t,x)=t^{-\alpha}V(xt^{-\beta})$ of the conservation law solves a first-order ODE, whose solutions are explicit up to a constant (the branches of Subsection~\ref{subsec:profile.eqs}). A profile is obtained by gluing branches along discontinuities that have to be admissible, that is, to satisfy the Rankine--Hugoniot and the Ole\u{\i}nik conditions (Subsection~\ref{subsect.adm}); in Subsection~\ref{subsect.construction} we choose the branches, locate the shocks and check the entropy conditions directly, and uniqueness follows from Theorem~\ref{unique.entropy}. Since $F$ changes convexity at the origin, the admissible jumps are not those of a convex flux, and their analysis leads to the threshold $\eps(\gamma,q)$.

\subsection{Organization of the paper}
Section~\ref{sect-Preliminaries} recalls the Marcinkiewicz spaces and the well-posedness theory, and Section~\ref{sect-estimates} contains the parabolic estimates, the one-sided estimates and the decay rates. Section~\ref{sect-CL} is devoted to uniqueness for problems with unbounded initial data. Theorems~\ref{farfield} and \ref{farfield.gamma1} are proved in Section~\ref{sect-farfield}, and Theorems~\ref{first.term.gamma<1} and \ref{first.term.gamma=1} in Sections~\ref{sect-lt1} and \ref{sect-eq1}, respectively. Section~\ref{sect-profiles} contains the proof of Theorem~\ref{self.profiles}, a complement describing the profiles of Theorem~\ref{first.term.gamma<1}\,(c) in detail.  The fixed-time compactness principle is proved in Appendix~\ref{sect-fixedtime}.

\section{Function spaces and well-posedness}\label{sect-Preliminaries}

Data satisfying \eqref{initial.gamma} are bounded but not integrable, and the natural space to measure them is the Marcinkiewicz space $L^{1/\gamma,\infty}(\rr)$. The first subsection recalls the Marcinkiewicz spaces $L^{p,\infty}(\rr)$ and the decreasing rearrangement, in terms of which the estimates of Section~\ref{sect-estimates} are most naturally expressed. The second one recalls the well-posedness theory for \eqref{eq.main} with bounded data in $L^{p_0}(\rr)$, which covers the data \eqref{initial.gamma} since they belong to $L^p(\rr)$ for every $p>\frac1\gamma$.

\subsection{Marcinkiewicz spaces and decreasing rearrangements}
For a measurable function $f$ on $\rr$ we denote by $d_f(\lambda)=|\{x\in\rr:|f(x)|>\lambda\}|$ its distribution function and by $f^*(s)=\inf\{\lambda>0:d_f(\lambda)\le s\}$ its decreasing rearrangement. For $1\le p<\infty$ the Marcinkiewicz space $L^{p,\infty}(\rr)$ consists of the functions $f$ such that
\begin{equation}\label{Lp.weak.norm}
\|f\|_{L^{p,\infty}(\rr)}:=\sup_{\lambda>0}\lambda\,d_f(\lambda)^{1/p}=\sup_{s>0}s^{1/p}f^*(s)<\infty .
\end{equation}

By \cite[Chapter~2, Theorem~6.2 and its proof]{BeSh88},
\begin{equation}\label{Kfunct}
    \int_0^sf^*(\sigma)\,d\sigma= \min _{{k\ge0}} \left(ks +\int_{\rr} {(|f|-k)_+}\,dx\right).
\end{equation}
with equality if $k=f^*(s)$, since then $f^*\ge k$ on $(0,s)$ and $f^*\le k$ on $(s,\infty)$.

On the other hand, if $E\subset\rr$ is measurable with $|E|=s<\infty$ (the case $|E|=\infty$ being trivial), then, for every $k\ge0$, $|f|\le k+(|f|-k)_+$ gives $\int_E|f|\,dx\le ks+\int_\rr(|f|-k)_+\,dx$. Choosing $k=f^*(s)$ and using the equality case in \eqref{Kfunct}, we obtain
\begin{equation}\label{HL}
\int_E|f|\,dx\le\int_0^{|E|}f^*(\sigma)\,d\sigma\qquad\text{for every measurable }E\subset\rr .
\end{equation}
For $1<p<\infty$, since $f^*$ is nonincreasing we have $\int_0^sf^*\ge s\,f^*(s)$, and, integrating $f^*(\sigma)\le\|f\|_{L^{p,\infty}(\rr)}\sigma^{-1/p}$, $\int_0^sf^*\le\frac p{p-1}\|f\|_{L^{p,\infty}(\rr)}\,s^{1-1/p}$. Hence
\begin{equation}\label{norm-equiv}
\|f\|_{L^{p,\infty}(\rr)}\le\sup_{s>0}s^{\frac1p-1}\int_0^sf^*(\sigma)\,d\sigma\le\frac p{p-1}\|f\|_{L^{p,\infty}(\rr)} .
\end{equation}
Finally, for $1\le p<r<\infty$ \cite[Proposition~1.1.14]{MR3243734}, 
\begin{equation}\label{interp.lorentz}
\|f\|_{L^r(\rr)}\le\Big(\frac r{r-p}\Big)^{1/r}\|f\|_{L^{p,\infty}(\rr)}^{p/r}\|f\|_{L^\infty(\rr)}^{1-p/r}.
\end{equation}

\subsection{Well-posedness}
The following result is classical. 

\begin{theorem}\label{integrable}
    Let $F\in C^1(\rr)$ with $F(0)=0$,  $1\le p_0<\infty$ and $u_0\in L^{p_0}(\rr)\cap L^\infty(\rr)$. Problem \eqref{eq.main} has a unique mild solution $u\in C([0,\infty);L^{p_0}(\rr))\cap L^\infty(Q)$, given by    \begin{equation}\label{duhamel}
        u(t)=G_t\ast u_0-\int_0^t\partial_xG_{t-s}\ast F(u(s))\,ds,\qquad t\ge0.
    \end{equation}
    Moreover:
    \begin{enumerate}
        \item[(i)] $u\in C((0,\infty);W^{2,p}(\rr))\cap C^1((0,\infty);L^p(\rr))$ for every $p\in(p_0,\infty)$; in particular $\partial_xu\in L^\infty((\tau,T)\times\rr)$ for $0<\tau<T$, and if $F'$ is locally H\"older continuous, $u$ is a classical solution in $Q$;
        \item[(ii)] for every convex $\eta\in C^1(\rr)$ with $\eta(0)=0$ and $|\eta(s)|\le C|s|^{p_0}$ for $|s|\le\|u_0\|_{L^\infty(\rr)}$, the function $t\mapsto\int_\rr\eta(u(t))\,dx$ is nonincreasing; in particular
            \[
                \|u(t)\|_{L^p(\rr)}\le\|u_0\|_{L^p(\rr)}\quad\textup{for }p\in[p_0,\infty];
            \]
        \item[(iii)] if $v_0\in L^{p_0}(\rr)\cap L^\infty(\rr)$ and $u_0\le v_0$ a.e., the corresponding solutions satisfy $u\le v$ in $Q$;
        \item[(iv)] if $u_{0,n}\to u_0$ in $L^{p_0}(\rr)$ and $\sup_n\|u_{0,n}\|_{L^\infty(\rr)}<\infty$, the corresponding solutions satisfy $u_n\to u$ in $C([0,T];L^{p_0}(\rr))$ for every $T>0$.
        \end{enumerate}
\end{theorem}

For $p_0=1$ Theorem~\ref{integrable} is proved in \cite[Theorem~5.1]{zuazuabook}, under the same assumptions on $F$; see also \cite[Proposition~1]{EZ} for $F(u)=|u|^{q-1}u$. For $p_0>1$ the proof is the same. Property (ii) follows by multiplying the equation in \eqref{eq.main} by $\eta'(u)$, after approximating $\eta$ by convex functions of class $C^2$, and observing that  $\int_\rr\eta'(u)\partial_x(F(u))\,dx=\int_\rr\partial_x\mathcal Q_\eta(u)\,dx=0$, where $\mathcal Q_\eta'=\eta'F'$ and $\mathcal Q_\eta(0)=0$. The bounds on $\|u(t)\|_{L^p(\rr)}$ follow by taking $\eta(s)=|s|^p$ if $p_0\le p<\infty$ and $p>1$, with $\eta_\delta(s)=(s^2+\delta^2)^{1/2}-\delta$ and $\delta\downarrow0$ if $p=p_0=1$, and by letting $p\to\infty$. Property (iv) follows from \eqref{duhamel} and Gronwall's lemma, and property (iii) from the case $p_0=1$ by approximation and (iv). If $F'$ is locally H\"older continuous, then, by (i), $F'(u)\partial_xu$ is locally H\"older continuous in $Q$, and $u$ is a classical solution by mollification and the interior Schauder estimate \cite[Theorem~4.9]{Lie96}. Classical regularity is used in the two-point maximum principle of Theorem~\ref{holder-oleinik} and in the comparison principles of Lemmas~\ref{ff.comparison} and \ref{n:comparison}.

\section{Parabolic and hyperbolic estimates}\label{sect-estimates}

This section collects the estimates that give the decay of the solutions. They are of two different natures. The first ones, in the next subsection, are of parabolic type, since they hold for every flux $F\in C^1(\rr)$ with $F(0)=0$, they come from the dissipation alone, and, for data in $L^{p_0,\infty}(\rr)$, they give the decay $t^{-\frac1{2p_0}}$ of the heat equation. The second ones are of hyperbolic type: they use the nonlinearity $F(u)=|u|^{q-1}u$, $q>1$, through one-sided (Ole\u{\i}nik-type) estimates, which limit how fast $u$ can increase to the right, and, combined with the $L^{p_0,\infty}$ bound, they give the rate $t^{-\frac1{p_0+q-1}}$ of the conservation law, which is better if and only if $q<p_0+1$. As explained in the Introduction, this competition is what decides the exponent $\theta$ in Theorem~\ref{first.term.gamma<1} and the critical curve $q=1+1/\gamma$.

\subsection{Parabolic estimates}

The parabolic estimates only use the dissipation of the equation, and they hold for every flux $F\in C^1(\rr)$ with $F(0)=0$. They are based on the following lemma, which controls the part of the solution above a level $k$ in terms of the same quantity for the datum; together with formula \eqref{Kfunct} for $\int_0^sf^*$, it will give in Theorem~\ref{week.integrable} that the solution does not increase the integral of the decreasing rearrangement of the datum. Level-set techniques for purely diffusive problems with Marcinkiewicz data were used in \cite{MR4671449}.

\begin{lemma}\label{levelk}
    Let $F\in C^1(\rr)$ with $F(0)=0$, let $u_0\in L^{p_0}(\rr)\cap L^\infty(\rr)$ for some $p_0\in[1,\infty)$, let $u$ be the solution of \eqref{eq.main} given by Theorem~\ref{integrable}, and let $k>0$. Then, for every $t>0$,
    \begin{align}
        \label{level-k}
        &\int_\rr(|u(t)|-k)_+^p\,dx\le\int_\rr(|u_0|-k)_+^p\,dx\qquad\text{for every }p\in[1,\infty),\\
        \label{level-k-infty}
        &\|(|u(t)|-k)_+\|_{L^\infty(\rr)}\le C_0\,t^{-1/2}\int_\rr(|u_0|-k)_+\,dx,
    \end{align}
    where $C_0$ is an absolute constant, independent of $F$.
\end{lemma}

\begin{proof}
All the integrals in \eqref{level-k} are finite, since $(|s|-k)_+^p\le k^{p-p_0}|s|^{p_0}$ if $p\le p_0$, and $(|s|-k)_+^p\le\|u_0\|_{L^\infty(\rr)}^{p-p_0}|s|^{p_0}$ for $|s|\le\|u_0\|_{L^\infty(\rr)}$ if $p\ge p_0$. Estimate \eqref{level-k} then follows from Theorem~\ref{integrable}(ii), applied for $p>1$ to the convex $C^1$ function $\eta(s)=(|s|-k)_+^p$, and for $p=1$ to the convex $C^1$ functions $\eta_\delta(s)=\big((|s|-k)_+^2+\delta^2\big)^{1/2}-\delta$, by monotone convergence as $\delta\downarrow0$.

Let $v:=(|u|-k)_+$. For $p\ge1$ the function $s\mapsto(|s|-k)_+^{2p}$ is convex and of class $C^1$, and the convective term in \eqref{eq.main} is a total derivative, so that
\[
    \frac{d}{dt}\int_\rr v^{2p}\,dx=-\frac{2(2p-1)}p\int_\rr\big|\partial_x\big(v^{p}\big)\big|^2dx,\qquad p\ge1.
\]
Moreover, $t\mapsto\int_\rr v^p(t)\,dx$ is nonincreasing, by \eqref{level-k}. We apply Nash's inequality in one dimension, $\|w\|_{L^2(\rr)}^6\le C\|w\|_{L^1(\rr)}^4\|\partial_xw\|_{L^2(\rr)}^2$, to $w=v^p$. Writing $y(t):=\int_\rr v^{2p}(t)\,dx$ and using that $\int_\rr v^p(t)\,dx\le\int_\rr v^p(s)\,dx$ for $t\ge s$, we get $y'\le-c_p\,y^3\big(\int_\rr v^p(s)\,dx\big)^{-4}$ for $t\ge s$, which integrates to
\[
    \|v(t)\|_{L^{2p}(\rr)}\le C_p\,(t-s)^{-\frac1{4p}}\,\|v(s)\|_{L^p(\rr)},\qquad 0\le s<t.
\]
Iterating this estimate for $p=2^n$, $n=0,1,\dots$, between the times $t(1-2^{-n})$ and $t(1-2^{-n-1})$, the powers of $t$ add up to $t^{-1/2}$ and the product of the remaining constants converges, which gives \eqref{level-k-infty}. See \cite[Theorem~5.2]{zuazuabook} (case $p=\infty$) for the details.
\end{proof}

\begin{theorem}\label{week.integrable}
    Let $F\in C^1(\rr)$ with $F(0)=0$,  $u_0\in L^\infty(\rr)\cap L^{p_0,\infty}(\rr)$ for some $p_0\in (1,\infty)$, and  $u$ be the solution of \eqref{eq.main} given by Theorem~\ref{integrable}.
    For every $t\ge0$ and $s>0$,
    \begin{equation}\label{HLP-dom}
        \int_0^su(t)^*(\sigma)\,d\sigma\le\int_0^su_0^*(\sigma)\,d\sigma .
    \end{equation}
    Consequently,
    \begin{align}
        \label{est.p.infty}
        &\|u(t)\|_{L^{p,\infty}(\rr)}\le\frac p{p-1}\,\|u_0\|_{L^{p,\infty}(\rr)}\qquad\text{for every }p\in[p_0,\infty),\ t\ge0,\\
        \label{decay.parabolic}
        &\|u(t)\|_{L^\infty(\rr)}\le C(p_0)\,\|u_0\|_{L^{p_0,\infty}(\rr)}\,t^{-\frac1{2p_0}}\qquad\text{for every }t>0.
    \end{align}

\end{theorem}

\begin{proof}
Let $s>0$ and $k>u_0^*(s)$. By \eqref{Kfunct} and \eqref{level-k},
\[
    \int_0^su(t)^*\,d\sigma\le ks+\int_\rr(|u(t)|-k)_+\,dx\le ks+\int_\rr(|u_0|-k)_+\,dx .
\]
Letting $k\downarrow u_0^*(s)$, by monotone convergence and the equality case of \eqref{Kfunct}, the right-hand side tends to $\int_0^su_0^*$, which proves \eqref{HLP-dom}. Multiplying \eqref{HLP-dom} by $s^{\frac1p-1}$ and taking the supremum in $s$, we get \eqref{est.p.infty} by \eqref{norm-equiv}.

Let $A:=\|u_0\|_{L^{p_0,\infty}(\rr)}$. For $k>0$,
\[
    \int_\rr(|u_0|-k)_+\,dx=\int_k^\infty d_{u_0}(\lambda)\,d\lambda\le A^{p_0}\int_k^\infty\lambda^{-p_0}\,d\lambda=\frac{A^{p_0}k^{1-p_0}}{p_0-1},
\]
so that, by \eqref{level-k-infty},
\[
    \|u(t)\|_{L^\infty(\rr)}\le k+\|(|u(t)|-k)_+\|_{L^\infty(\rr)}\le k+\frac{C_0A^{p_0}}{(p_0-1)t^{1/2}}\,k^{1-p_0}.
\]
The choice $k=C_0^{1/p_0}At^{-1/(2p_0)}$ gives \eqref{decay.parabolic}.
\end{proof}

\begin{remark}\label{rem.marc}
    (a) By \eqref{HL}, $\int_0^sf^*$ is the largest integral of $|f|$ over a set of measure $s$, hence \eqref{HLP-dom} says that, for every $s>0$, the integral of $|u(t)|$ over any set of measure $s$ is at most the largest integral of $|u_0|$ over a set of measure $s$. 

    \noindent (b) The space $L^{1,\infty}(\rr)$ is not preserved by the heat semigroup, not even for bounded data. One can construct $u_0\in L^{1,\infty}(\rr)\cap L^\infty(\rr)$ such that $\|G_t\ast u_0\|_{L^{1,\infty}(\rr)}=\infty$ for every $t>0$. For such a datum the solution of \eqref{eq.main} also satisfies $\|u(t)\|_{L^{1,\infty}(\rr)}=\infty$ for every $t>0$.

However, for data as in \eqref{initial.gamma} with $\gamma=1$, inequality \eqref{HLP-dom} still holds.
It yields the logarithmic bound \eqref{local.L1.est} of Corollary~\ref{local.L1} below. This, together with Proposition~\ref{cor.decay.gamma} and the far-field analysis in Section~\ref{sect-farfield}, will show that  $u(t)\in L^{1,\infty}(\rr)$ for all times.

\end{remark}

For data as in \eqref{initial.gamma}, \eqref{HL} and \eqref{HLP-dom} give local $L^1$ bounds that are uniform in time.

\begin{cor}\label{local.L1}
    Let $F\in C^1(\rr)$ with $F(0)=0$, let $\gamma\in(0,1]$, let $u_0\in L^\infty(\rr)$ satisfy $|u_0(x)|\le\kappa(1+|x|)^{-\gamma}$ for a.e.\ $x\in\rr$, and let $u$ be the solution of \eqref{eq.main} given by Theorem~\ref{integrable}. Then $u_0\in L^{p,\infty}(\rr)$ for every $p\ge\frac1\gamma$, with $\|u_0\|_{L^{1/\gamma,\infty}(\rr)}\le2^\gamma\kappa$, and for every measurable set $E\subset\rr$ and every $t\ge0$,
    \begin{equation}\label{local.L1.est}
        \int_E|u(t,x)|\,dx\le\int_0^{|E|}u_0^*(\sigma)\,d\sigma\le
        \begin{cases}
            \dfrac{2^\gamma\kappa}{1-\gamma}\,|E|^{1-\gamma},&0<\gamma<1,\\[8pt]
            2\kappa\log\big(1+\tfrac{|E|}2\big),&\gamma=1.
        \end{cases}
    \end{equation}
\end{cor}

\begin{proof}The first inequality in \eqref{local.L1.est} follows from \eqref{HL} and \eqref{HLP-dom}.
For such initial data  we have $u_0^*(\sigma)\le\kappa(1+\sigma/2)^{-\gamma}$. Integrating this inequality gives the desired estimate.
\end{proof}

\subsection{One-sided estimates} We recall and extend the one-sided estimates known as Ole\u{\i}nik-type estimates \cite{Ole63}.

\begin{proposition}[\cite{Ole63,EVZArma}]\label{oleinik.prop}
    Let $1<q\le2$ and  $u\ge0$ be a bounded classical solution of the equation in \eqref{eq.main}. Then
    \begin{equation}\label{oleinik}
        u^{q-1}(t,x)-u^{q-1}(t,y)\le\frac{x-y}{qt}\qquad\text{for all }t>0,\ y<x.
    \end{equation}
\end{proposition}

For $1<q<2$  estimate \eqref{oleinik} is proved in \cite[Lemma~1.1]{EVZArma}. For $q=2$ it is Ole\u{\i}nik's classical one-sided estimate \cite{Ole63}; see also \cite[Proposition~4]{watanabe}. 

  For $q>2$ we replace the above estimate  by the following one, which coincides with \eqref{oleinik} when $q=2$.

\begin{theorem}\label{holder-oleinik}
    Let $q\ge2$ and let $u\ge0$ be a bounded classical solution of the equation in \eqref{eq.main} such that $\partial_xu\in L^\infty((\tau,T)\times\rr)$ for all $0<\tau<T$. Then
    \begin{equation}\label{holder-oleinik-est}
        u(t,x)-u(t,y)\le\Big(\frac{x-y}{qt}\Big)^{\frac1{q-1}}\qquad\text{for all }t>0,\ y<x .
    \end{equation}
\end{theorem}

The proof relies on a two-point maximum principle. Two-point functions have been used in the analysis of parabolic equations in one space variable since the work of Kru\v{z}kov \cite{Kru67}; see \cite{MR2517770,MR2504416,MR3381494,MR4294711} and the references therein.

The idea is to compare $u(t,x)-u(t,y)$ with a function of $t$ and $x-y$ that blows up at an earlier time. At a positive maximum of the difference the first-order conditions give the same slope of $u$ at $x$ and at $y$, the diffusion has a favourable sign, and the convection produces a negative term because of the inequality 
\begin{equation}
  \label{eq:ineg.tonta}
  a^{q-1}-b^{q-1}\ge(a-b)^{q-1}, \quad a\ge b\ge0, \ q\geq 2.
\end{equation}
This is the only place where the condition $q\ge2$ is used.

\begin{proof}
 Let $K>q^{-\frac1{q-1}}$ and set
\[
    \omega(t,z):=K\Big(\frac z{t-s}\Big)^{\frac1{q-1}},\qquad 0<s<t\leq T,\ z>0.
\]
It is sufficient to prove that
\begin{equation}\label{holder.delayed}
    S:=\sup_{s<t\le T,\ y<x}\big(u(t,x)-u(t,y)-\omega(t,x-y)\big)\le0.
\end{equation}
Indeed, since $T$ is arbitrary, \eqref{holder-oleinik-est} follows from \eqref{holder.delayed} by letting $s\downarrow0$ and then $K\downarrow q^{-\frac1{q-1}}$.

Assume by contradiction that $S>0$. Since $u\ge0$, we have $S\le\|u\|_{L^\infty(Q)}$. Let $L:=\|\partial_xu\|_{L^\infty((s,T)\times\rr)}$.
We choose a point $(t_*,x_*,y_*)$, with $s<t_*\le T$ and $y_*<x_*$, at which
\[u(t_*,x_*)-u(t_*,y_*)-\omega(t_*,x_*-y_*)>\frac{3S}4.\]
 For $\psi(x):=(1+x^2)^{1/2}$ we choose $\eps$, which will be fixed later, such that
 \[\eps(\psi(x_*)+\psi(y_*))<\frac S4,\]
  and we introduce
\[
    \Phi(t,x,y):=u(t,x)-u(t,y)-\omega(t,x-y)-\eps\big(\psi(x)+\psi(y)\big),
\]
so that $\Phi(t_*,x_*,y_*)>\frac S2$.

We now prove that $\Phi$ attains its supremum away from the diagonal. Let $s<t\le T$ and $y<x$. Since $\omega\ge0$ and $\psi\ge0$,
\[
    \Phi(t,x,y)\le u(t,x)-u(t,y)\le L(x-y)\le\frac S2\qquad\text{if }x-y\le\frac S{2L}.
\]
Since $\omega$ is increasing in $z$ and $u\ge0$,
\[
    \Phi(t,x,y)\le\|u\|_{L^\infty(Q)}-\omega\Big(t,\frac S{2L}\Big)-\eps\big(\psi(x)+\psi(y)\big)\qquad\text{if }x-y\ge\frac S{2L}.
\]
The right-hand side tends to $-\infty$ as $t\downarrow s$ and as $|x|+|y|\to\infty$. Hence there are $\rho,R>0$ such that $\{\Phi>\frac S2\}$ is contained in the compact set
\[
    \mathcal K:=\Big\{s+\rho\le t\le T,\ |x|+|y|\le R,\ x-y\ge\frac S{2L}\Big\}.
\]
Thus $\Phi$ attains its positive supremum over $\{s<t\le T,\ y<x\}$ at a point $(t_0,x_0,y_0)\in\mathcal K$. We denote
\[
    z_0:=x_0-y_0\ge\frac S{2L},\qquad\omega_0:=\omega(t_0,z_0),\qquad a:=u(t_0,x_0),\qquad b:=u(t_0,y_0).
\]
Since $\Phi(t_0,x_0,y_0)>0$, $\psi\ge0$ and $u\ge0$, we have
\begin{equation}\label{holder.gap}
    a-b>\omega_0>0,\qquad b\ge0.
\end{equation}

Denoting $\partial_z\omega_0:=\partial_z\omega(t_0,z_0)$, the first-order conditions in $x$ and $y$ read
\[
    \partial_xu(t_0,x_0)=\partial_z\omega_0+\eps\psi'(x_0),\qquad \partial_xu(t_0,y_0)=\partial_z\omega_0-\eps\psi'(y_0).
\]
The map $r\mapsto\Phi(t_0,x_0+r,y_0+r)$ attains its maximum at $r=0$, hence the second derivative of this function at $r=0$ is nonpositive. Since $0\le\psi''\le1$, we get
\[
    \partial_{xx}u(t_0,x_0)-\partial_{xx}u(t_0,y_0)\le\eps\big(\psi''(x_0)+\psi''(y_0)\big)\le2\eps.
\]
Since $s<t_0\le T$, we also have $\partial_t\Phi(t_0,x_0,y_0)\ge0$, that is,
\[
 0\leq \partial_t\Phi (t_0,x_0,y_0)=   \partial_tu(t_0,x_0)-\partial_tu(t_0,y_0)-\partial_t \omega(t_0,z_0).
\]
We insert in this inequality the equation $\partial_tu=\partial_{xx}u-qu^{q-1}\partial_xu$. Using the above relations on the first-order derivatives, $|\psi'|\le1$ and $0\le b<a\le\|u\|_{L^\infty(Q)}$, we obtain, using the explicit expression of $\omega$, that
\begin{align*}
    0&\le\partial_{xx}u(t_0,x_0)-\partial_{xx}u(t_0,y_0)-qa^{q-1}\partial_xu(t_0,x_0)\\
&\quad+qb^{q-1}\partial_xu(t_0,y_0)+\frac{\omega_0}{(q-1)(t_0-s)}\\
&\le 2 \eps - qa^{q-1}
(\partial_z\omega_0+\eps\psi'(x_0))+qb^{q-1}(\partial_z \omega_0-\eps \psi'(y_0))+\frac{\omega_0}{(q-1)(t_0-s)}\\
    &\le \eps(2+2q\|u\|_{L^\infty(Q)}^{q-1})-q\big(a^{q-1}-b^{q-1}\big)\partial_z\omega_0+\frac{\omega_0}{(q-1)(t_0-s)}\\
    &\le \eps(2+2q\|u\|_{L^\infty(Q)}^{q-1})-q\omega_0^{q-1}\partial_z\omega_0+\frac{\omega_0}{(q-1)(t_0-s)},
\end{align*}
where we used that, for
 $q\ge2$, \eqref{holder.gap} and \eqref{eq:ineg.tonta} yield $ 
    a^{q-1}-b^{q-1}\ge(a-b)^{q-1}>\omega_0^{\,q-1}.
$
Therefore, using that $\omega_0^{q-2}\partial_z\omega_0 =\frac{K^{q-1}}{(q-1)(t_0-s)}$, $z_0\ge\frac S{2L}$ and $t_0-s\le T-s$, we obtain \begin{align*}
    0&\le \eps(2+2q\|u\|_{L^\infty(Q)}^{q-1})-\big(qK^{q-1}-1\big)\frac{\omega_0}{(q-1)(t_0-s)}\\
    &=\eps(2+2q\|u\|_{L^\infty(Q)}^{q-1})-\frac{K\big(qK^{q-1}-1\big)}{q-1}\,z_0^{\frac1{q-1}}(t_0-s)^{-\frac q{q-1}}\\
   &\leq  \eps(2+2q\|u\|_{L^\infty(Q)}^{q-1})-\frac{K\big(qK^{q-1}-1\big)}{q-1}\,\big(\frac {S}{2L}\big)^{\frac1{q-1}}(T-s)^{-\frac q{q-1}} <0,
\end{align*}
which gives us a contradiction by choosing $\eps$ small enough. Hence $S\le0$, and this finishes the proof of Theorem~\ref{holder-oleinik}.
\end{proof}

\begin{remark}\label{rem.hoff}
Estimate \eqref{holder-oleinik-est} is invariant under the scaling $(t,x)\mapsto(t/\nu,x/\nu)$, which transforms \eqref{eq.main} into $u_t+(F(u))_x=\nu u_{xx}$; it is therefore uniform with respect to the viscosity. Let $u$ be a nonnegative entropy solutions of the conservation law \eqref{conservation.law}. The convexity of $F$ on $[0,\infty)$ and the sharp Ole\u{\i}nik estimate of Hoff \cite[Theorem~1]{Hof83} in \eqref{oleinik}, together with the inequality \eqref{eq:ineg.tonta}
 give \eqref{holder-oleinik-est} with the same constant for every $q\ge2$.
\end{remark}

\subsection{From one-sided estimates to decay}

The one-sided estimates are turned into decay by the following elementary lemma. It says that a function in $L^{p,\infty}(\rr)$ that cannot increase too fast to the right is bounded, with a bound that only depends on its $L^{p,\infty}$ norm and on the constant in the one-sided estimate. Applied to $u(t)$, it gives the decay of the solution in the convective regime (Theorem~\ref{th.decay.p_0}).

\begin{lemma}\label{n:magic3}
    Let $p\in[1,\infty)$, $r,\sigma>0$, $A,B>0$, and let $u$ be a measurable function on $\rr$ with $\|u\|_{L^{p,\infty}(\rr)}\le A$ and
    \begin{equation}\label{H2}
        \big(|u|^{r-1}u\big)(x)-\big(|u|^{r-1}u\big)(y)\le\big(B(x-y)\big)^\sigma\qquad\text{for a.e.\ }(x,y)\in\rr^2\text{ with }y<x.
    \end{equation}
    Then
    \begin{equation}\label{n:magic3-est}
        \|u\|_{L^\infty(\rr)}\le C(p,r,\sigma)\,A^{\frac{p\sigma}{p\sigma+r}}B^{\frac\sigma{p\sigma+r}}.
    \end{equation}
\end{lemma}

\begin{proof}
Let $0<m<\|u\|_{L^\infty(\rr)}$, and let $x_0$ be such that $|u(x_0)|>m$ and \eqref{H2} holds with $x=x_0$ for a.e.\ $y<x_0$ and with $y=x_0$ for a.e.\ $x>x_0$. Fix $\mu\in(0,1)$ and let
\[
    \rho:=\frac1B\big((1-\mu^r)\,m^r\big)^{\frac1\sigma},
\]
so that $(B\rho)^\sigma=(1-\mu^r)\,m^r$. If $u(x_0)>m$, then, by \eqref{H2}, for a.e.\ $y\in(x_0-\rho,x_0)$,
\[
    \big(|u|^{r-1}u\big)(y)\ge u^r(x_0)-\big(B(x_0-y)\big)^\sigma>m^r-(B\rho)^\sigma=\mu^rm^r,
\]
so that $u>\mu m$ a.e.\ in $(x_0-\rho,x_0)$. If $u(x_0)<-m$, then, by \eqref{H2}, for a.e.\ $x\in(x_0,x_0+\rho)$,
\[
    \big(|u|^{r-1}u\big)(x)\le-|u(x_0)|^r+\big(B(x-x_0)\big)^\sigma<-m^r+(B\rho)^\sigma=-\mu^rm^r,
\]
so that $u<-\mu m$ a.e.\ in $(x_0,x_0+\rho)$. In both cases $d_u(\mu m)\ge\rho$ (see \eqref{Lp.weak.norm}). Hence, by \eqref{Lp.weak.norm},
\[
    A\ge\mu m\,d_u(\mu m)^{1/p}\ge\mu m\,\rho^{1/p}=\mu(1-\mu^r)^{\frac1{p\sigma}}\,m^{\frac{p\sigma+r}{p\sigma}}B^{-\frac1p}.
\]
Choosing $\mu=\big(\frac{p\sigma}{p\sigma+r}\big)^{1/r}$, which maximizes $\mu(1-\mu^r)^{1/(p\sigma)}$, and letting $m\uparrow\|u\|_{L^\infty(\rr)}$, we obtain \eqref{n:magic3-est}.
\end{proof}

For initial data $u_0\in L^1(\rr)\cap L^\infty(\rr)$,  problem \eqref{eq.main} has a unique solution \cite[Theorem~5.1, p.~53]{zuazuabook} in $C([0,\infty),L^1(\rr))\cap L^\infty((0,\infty)\times \rr)$, which satisfies
\[
\|u(t)\|_{L^{1,\infty}(\rr)}\leq \|u(t)\|_{L^{1}(\rr)}\leq \|u_0\|_{L^1(\rr)}.
\]
In this case, Proposition \ref{oleinik.prop}, Theorem \ref{holder-oleinik} and Lemma~\ref{n:magic3} provide a hyperbolic estimate that holds for all $q > 1$:
\[
\|u(t)\|_{L^\infty(\rr)}\lesssim \|u_0\|_{L^1(\rr)}^{\frac{1}{q}}\,t^{-\frac1{q}} \quad \text{for all}\, t>0.
 \]
 This allows one to extend the   $L^1(\rr)$ well-posedness theory beyond $q=2$, see
 \cite[Remark~5.1 and Theorem~5.3]{zuazuabook}.

\begin{theorem}\label{th.decay.p_0}
    Let $p_0\in(1,\infty)$, $q>1$, $u_0\in L^{p_0,\infty}(\rr)\cap L^\infty(\rr)$, $A:=\|u_0\|_{L^{p_0,\infty}(\rr)}$, and let $u$ be the solution of \eqref{eq.main} given by Theorem~\ref{integrable}. Then, for every $t>0$,
    \begin{equation}\label{decay}
        \|u(t)\|_{L^\infty(\rr)}\le\min\Big\{C_1(p_0)\,A\,t^{-\frac1{2p_0}},\ C_2(p_0,q)\,A^{\frac{p_0}{p_0+q-1}}\,t^{-\frac1{p_0+q-1}}, \|u_0\|_{L^\infty(\rr)}\Big\}.
    \end{equation}

\end{theorem}

The first bound is the parabolic one and does not see the nonlinearity; the second one is of hyperbolic nature, and it gives the better rate for large times if and only if $q<p_0+1$. The hyperbolic rate is that of the inviscid problem; for convex fluxes, such rates were derived by Dafermos for entropy solutions of conservation laws from the growth at infinity of the primitive of the datum \cite[Theorem~11.5.1]{Daf2016}, \cite[Theorem~7.1]{Daf77}.

 For $F(u)=|u|^{q-1}u$ the well-posedness theory for \eqref{eq.main} can be developed for  $u_0$ in $L^{p_0,\infty}(\rr)$, $p_0>1$, not necessarily bounded. This can be done using the hyperbolic estimates for $t$ close to zero.

\begin{proof}
The first bound is \eqref{decay.parabolic}. For the second one, let $w$ be the solution of \eqref{eq.main} with initial datum $|u_0|$.  Theorem~\ref{integrable}(iii) gives $|u|\le w$. By \eqref{est.p.infty}, $\|w(t)\|_{L^{p_0,\infty}(\rr)}\le\frac{p_0}{p_0-1}A$. If $1<q\le2$, Proposition~\ref{oleinik.prop} and Lemma~\ref{n:magic3} with $(r,\sigma)=(q-1,1)$ and $B=\frac1{qt}$ give the bound for $w(t)$, hence for $u(t)$. If $q\ge2$, Theorem~\ref{holder-oleinik}, which applies to $w$ by Theorem~\ref{integrable}(i), and Lemma~\ref{n:magic3} with $(r,\sigma)=(1,\frac1{q-1})$ and $B=\frac1{qt}$ give the same conclusion. In both cases $\frac{p_0\sigma}{p_0\sigma+r}=\frac{p_0}{p_0+q-1}$ and $\frac\sigma{p_0\sigma+r}=\frac1{p_0+q-1}$.
The $L^\infty(\rr)$ bound is given by Theorem~\ref{integrable}(ii).
\end{proof}

 For data satisfying \eqref{initial.gamma} with  $0<\gamma<1$,  we can apply Theorem~\ref{th.decay.p_0}
with $p_0=1/\gamma$ to obtain
the decay rate stated in Theorem~\ref{first.term.gamma<1}.
For $\gamma=1$ we can also obtain
 the decay rates stated in Theorem~\ref{first.term.gamma=1} by using the results in Corollary~\ref{local.L1}, as we prove next.

\begin{proposition}
\label{cor.decay.gamma}
    Let $q>1$, $\gamma=1$, let $u_0$ satisfy \eqref{initial.gamma}
     and let \[\kappa:=\esssup_{x\in\rr}(1+|x|)|u_0(x)|.\]  Then \eqref{thm2.decay} holds, with constants depending only on $q$ and $\kappa$.
\end{proposition}

\begin{proof}
Since $\|u_0\|_{L^\infty(\rr)}\le\kappa$, Theorem~\ref{integrable}(ii) gives $\|u(t)\|_{L^\infty(\rr)}\le\kappa$ for all $t\ge0$, which is enough for $t\le1$. Let $t\ge1$.

 If $q>2$, we use \eqref{level-k-infty}. Since $|u_0(x)|\le\kappa(1+|x|)^{-1}$, for $0<k\le\kappa$ we have
\[
\int_\rr(|u_0|-k)_+\,dx\le2\int_0^{\frac\kappa k-1}\Big(\frac\kappa{1+x}-k\Big)\,dx\le2\kappa\log\frac\kappa k ,
\]
and the choice $k=\kappa t^{-1/2}$ gives
\[
\|u(t)\|_{L^\infty(\rr)}\le k+\|(|u(t)|-k)_+\|_{L^\infty(\rr)}\le\kappa t^{-1/2}\big(1+C_0\log t\big),
\]
which is \eqref{thm2.decay}.

If $1<q\le2$, let $w$ be the solution of \eqref{eq.main} with initial datum $|u_0|$. Then $|u|\le w$ by Theorem~\ref{integrable}(iii), and $\|w(t)\|_{L^\infty(\rr)}\le\kappa$ by Theorem~\ref{integrable}(ii). Let $0<m<\|w(t)\|_{L^\infty(\rr)}$, let $x_0\in\rr$ be such that $w(t,x_0)>m$, and let
\[
c:=1-2^{1-q},\qquad \rho:=cqt\,m^{q-1},\qquad I:=(x_0-\rho,x_0).
\]
By \eqref{oleinik}, for every $y\in I$,
\[
w^{q-1}(t,y)\ge w^{q-1}(t,x_0)-\frac{x_0-y}{qt}>m^{q-1}-\frac\rho{qt}=\Big(\frac m2\Big)^{q-1}.
\]
Hence $w(t)>\frac m2$ on $I$, and \eqref{local.L1.est}, written for $w$ with $E=I$, gives
\[
\frac m2\,\rho\le\int_Iw(t,x)\,dx\le2\kappa\log\Big(1+\frac\rho2\Big).
\]
Since $m<\kappa$ and $q\le2$, we have $\frac\rho2\le t\kappa^{q-1}$. Therefore
\[
\frac{cq}2\,t\,m^q\le2\kappa\log\big(1+t\kappa^{q-1}\big)\le C\log(2+t),
\]
and letting $m\uparrow\|w(t)\|_{L^\infty(\rr)}$ we obtain \eqref{thm2.decay}, since $|u|\le w$.
\end{proof}

\section{Uniqueness for unbounded initial data}\label{sect-CL}

In this section we prove uniqueness for the conservation law \eqref{conservation.law} and for the viscous equation \eqref{eq.main} with locally integrable data that may be unbounded, as $U_0$ is, in classes defined by a rate of blow-up of $\|u(t)\|_{L^\infty(\rr)}$ as $t\downarrow0$ (Theorem~\ref{unique.entropy}). This is what identifies the limits of the rescaled solutions in Section~\ref{sect-lt1}  and the profiles of Theorem~\ref{self.profiles}. Since the data are unbounded, the classical theory for bounded data does not apply.

The results of this section hold for fluxes more general than powers. We assume that $F:\rr\to\rr$ is a locally Lipschitz function such that $F(0)=0$ and, for some $q>1$ and $C_F>0$,
\begin{equation}\label{F.growth}
|F(u)-F(v)|\le C_F\big(|u|^{q-1}+|v|^{q-1}\big)|u-v|\quad\text{for all }u,v\in\rr ;
\end{equation}
in particular $|F(u)|\le C_F|u|^q$. The flux $F(u)=|u|^{q-1}u$ satisfies \eqref{F.growth} with $C_F=q$. Equations~\eqref{eq.main} and \eqref{conservation.law} are considered here for this general flux $F$.

We now define solutions for initial data that are only in $L^1_{\textup{loc}}(\rr)$. Since solutions may be unbounded as $t\downarrow0$, we only require them to be bounded away from $t=0$, and we ask them to take the datum in $L^1_{\textup{loc}}(\rr)$.

\begin{definition}\label{def.entropy}
    Let $u_0\in L^1_{\textup{loc}}(\rr)$.
    \begin{enumerate}
        \item[(a)] A function $u$ is an \emph{entropy solution} of \eqref{conservation.law} if $u\in L^\infty((\tau,\infty)\times\rr)$ for every $\tau>0$,
        \begin{equation}
            \label{ineg.k}
            \partial_t|u-k|+\partial_x\big(\sgn(u-k)(F(u)-F(k))\big)\le0\quad\text{in }\mathcal D'(Q)\text{ for every }k\in\rr,
            \end{equation}
            and
            \begin{equation}
          \label{initial.entropy}    \esslim_{t\downarrow0}\int_K|u(t,x)-u_0(x)|\,dx=0\quad\text{for every compact set }K\subset\rr .
        \end{equation}
        \item[(b)] A function $u$ is a \emph{very weak solution} of \eqref{eq.main} if $u\in L^\infty((\tau,\infty)\times\rr)$ for every $\tau>0$, the equation holds in $\mathcal D'(Q)$, and \eqref{initial.entropy} holds.
    \end{enumerate}
\end{definition}

For bounded data, existence and uniqueness of bounded entropy solutions were established by Kru\v{z}kov \cite{MR267257}. For unbounded data, such as $U_0$,  one could work with the notion of renormalized solution introduced by B\'enilan, Carrillo and Wittbold \cite{MR1784177}; see also \cite{MR2839869}.
  However, determining if our limit solutions are renormalized or not may be difficult. Hence we prefer instead to impose  a rate of blow-up of $\|u(t)\|_{L^\infty(\rr)}$ as $t\downarrow0$ in order to obtain uniqueness. In the proof we
 shall use the following classical result for the difference of two solutions.
\begin{lemma}\label{u:kato}
    Let $0<\tau<T$, $\nu\ge0$. Let    $u_1,u_2\in L^\infty((\tau,T)\times\rr)$
     satisfy 
    \begin{equation}\label{visc.eq.nu}
        u_t+(F(u))_x=\nu u_{xx}\quad\text{in }\mathcal D'((\tau,T)\times\rr)
    \end{equation}
    for some $\nu>0 $, or ($\nu=0$)
         \[\partial_t|u-k|+\partial_x\big(\sgn(u-k)(F(u)-F(k))\big)\le0\quad\text{in }\mathcal D'((\tau,T)\times\rr)\text{ for every }k\in\rr.\]

    Then, in $\mathcal D'((\tau,T)\times\rr)$,
    \begin{equation}\label{kato.ineq}
        \partial_t|u_1-u_2|+\partial_x\big(\sgn(u_1-u_2)(F(u_1)-F(u_2))\big)\le\nu\,\partial_{xx}|u_1-u_2| .
    \end{equation}
    Consequently, for every $\phi\in C^2_c(\rr)$ with $\phi\ge0$, and for a.e.\ $\tau<t_1<t_2<T$ \textup{(}for all such $t_1,t_2$ if $\nu>0$\textup{)},
    \begin{equation}\label{entropy.for.two}
        \begin{aligned}
            \int_\rr|u_1-u_2|(t_2)\,\phi\,dx&\le\int_\rr|u_1-u_2|(t_1)\,\phi\,dx\\
            &\quad+\int_{t_1}^{t_2}\!\!\int_\rr\big(|F(u_1)-F(u_2)|\,|\phi'|+\nu|u_1-u_2|\,|\phi''|\big)\,dx\,ds.
        \end{aligned}
    \end{equation}
\end{lemma}

\begin{proof}
For $\nu=0$, \eqref{kato.ineq} is Kru\v{z}kov's inequality, obtained by the doubling of variables \cite[inequality~(3.7)]{MR267257}; since $u_1$ and $u_2$ are bounded on the slab, only the Lipschitz constant of $F$ on a bounded interval is involved. For $\nu>0$ the inequality is classical. Since $F(u_i)$ is bounded, the local energy estimate gives $\partial_xu_i\in L^2_{\textup{loc}}((\tau,T)\times\rr)$, and one tests the equation satisfied by $u_1-u_2$ with $S_\delta'(u_1-u_2)\Phi$, where $\Phi\ge0$ is a test function and $S_\delta$ are convex approximations of the absolute value, and lets $\delta\downarrow0$; see, e.g., \cite[Proposition~4.2]{endal}.  Finally,~\eqref{entropy.for.two} follows from~\eqref{kato.ineq} by testing with $\phi(x)\psi(t)$, where $\psi\ge0$ approximates the indicator function of $[t_1,t_2]$, at the Lebesgue points $t_1$, $t_2$ of $t\mapsto\int_\rr|u_1(t)-u_2(t)|\phi\,dx$, and integrating by parts in $x$; for $\nu>0$ the solutions are continuous in time with values in $L^1_{\textup{loc}}(\rr)$, so that every $t_1<t_2$ is admissible.
\end{proof}

\begin{remark}\label{kato.function}
    Let $\nu>0$,  $u\in L^\infty((\tau,T)\times\rr)$ satisfy \eqref{visc.eq.nu}, and  let $\varsigma\in C^2(\rr)$ be bounded together with its first two derivatives. Then $z:=u-\varsigma$ satisfies
    \[
        \partial_tz+\partial_x\big(F(u)-F(\varsigma)\big)-\nu\,\partial_{xx}z=-(F(\varsigma))'+\nu\varsigma''\quad\text{in }\mathcal D'((\tau,T)\times\rr),
    \]
    and the same argument above, applied to $z$, gives \eqref{entropy.for.two} with $u_1=u$ and $u_2=\varsigma$, with the additional term $\int_{t_1}^{t_2}\!\int_\rr\big(|(F(\varsigma))'|+\nu|\varsigma''|\big)\phi\,dx\,ds$ on the right-hand side.
\end{remark}

The proof of   uniqueness  is not difficult. The entropy inequality for the difference of two solutions (Lemma~\ref{u:kato}) gives a recursive inequality between the $L^1$ norms of the difference on balls of radii $R$ and $2R$, and iterating it we show that these norms vanish. The same argument covers the conservation law and the viscous equation, and every flux with the growth \eqref{F.growth}.

\begin{theorem}\label{unique.entropy}
    Let $q>1$,  $F$ satisfying \eqref{F.growth},  $u_0\in L^1_{\textup{loc}}(\rr)$, and  $0<\alpha<\frac1q$. Then there exists at most one entropy solution $u$ of \eqref{conservation.law}, and at most one very weak solution $u$ of \eqref{eq.main}, such that
    \begin{equation}\label{imposed.decay.old}
        \esslimsup_{t\downarrow0}\ t^{\alpha}\|u(t)\|_{L^\infty(\rr)}<\infty.
    \end{equation}
\end{theorem}

\begin{remark}\label{rem.uniq}
    For data in $L^{p_0,\infty}(\rr)$, $p_0>1$, the natural exponent is $\alpha=\frac1{p_0+q-1}<\frac1q$, the rate of the hyperbolic bound in \eqref{decay}.
\end{remark}

\begin{proof}
Let $u_1$ and $u_2$ be two solutions, both entropy solutions of \eqref{conservation.law} or both very weak solutions of~\eqref{eq.main}, satisfying \eqref{imposed.decay.old}. We set $\nu=0$ in the first case and $\nu=1$ in the second one, and fix $T>0$. By~\eqref{imposed.decay.old}, and since $u_i\in L^\infty((\tau,\infty)\times\rr)$ for every $\tau>0$, there is a constant $C_T$ such that
\begin{equation}\label{imposed.decay}
    |u_i(t,x)|\le C_Tt^{-\alpha}\quad\text{for a.e.\ }(t,x)\in(0,T)\times\rr,\ i=1,2,
\end{equation}
and then, by \eqref{F.growth},
\[
    |F(u_1)-F(u_2)|\le2C_FC_T^{q-1}t^{-\alpha(q-1)}|u_1-u_2|\quad\text{a.e. in}\ (0,T)\times\rr.
\]

Let $\varphi\in C^\infty_c(\rr)$ with $0\le\varphi\le1$, $\varphi=1$ on $[-1,1]$, $\operatorname{supp}\varphi\subset[-2,2]$, $|\varphi'|\le2$ and $|\varphi''|\le8$, let $\varphi_R:=\varphi(\cdot/R)$, and let \[D_R(t):=\int_{|x|<R}|u_1(t)-u_2(t)|\,dx.\]
 By \eqref{entropy.for.two} with $\phi=\varphi_R$, for $R\ge1$ and a.e.\ $t_1<t$ in $(0,T)$,
\[
    D_R(t)\le\int_\rr|u_1(t_1)-u_2(t_1)|\varphi_R\,dx+\int_{t_1}^tk_R(s)\,D_{2R}(s)\,ds,\qquad k_R(s):=\frac{4C_FC_T^{q-1}}{R\,s^{\alpha(q-1)}}+\frac{8\nu}{R^2}.
\]
Letting $t_1\downarrow0$ along a sequence for which \eqref{initial.entropy} holds for $u_1$ and $u_2$, the first term on the right-hand side tends to $0$, since it is bounded by $\sum_{i=1,2}\int_{|x|<2R}|u_i(t_1)-u_0|\,dx$,
and, by monotone convergence,
\begin{equation}\label{u:rec}
    D_R(t)\le\int_0^tk_R(s)\,D_{2R}(s)\,ds\qquad\text{for a.e.\ }t\in(0,T)\text{ and all }R\ge1 .
\end{equation}

Let $\mathcal D(R):=\esssup_{t\in(0,T)}t^{\alpha} D_R(t)$. By \eqref{imposed.decay}, $D_R(t)\le4 C_TR t^{-\alpha}$, so that $\mathcal D(R)\le4C_TR$. Since $D_{2R}(s)\le s^{-\alpha}\mathcal D(2R)$ and $\alpha q<1$, \eqref{u:rec} gives
\[
    \mathcal D(R)\le\frac{K_T}R\,\mathcal D(2R)\qquad\text{for }R\ge1,\qquad K_T:=\frac{4C_FC_T^{q-1}}{1-\alpha q}\,T^{1-\alpha(q-1)}+\frac{8\nu T}{1-\alpha}.
\]
Applying this inequality at the radii $R,2R,\dots,2^{n-1}R$, we obtain
\[
    \mathcal D(R)\le\frac{K_T^n}{R^n\,2^{n(n-1)/2}}\,\mathcal D(2^nR)\le4C_T\,K_T^n\,2^{n-\frac{n(n-1)}2}\,R^{1-n}\xrightarrow[n\to\infty]{}0 .
\]
Hence $\mathcal D\equiv0$ on $[1,\infty)$, that is, $u_1=u_2$ a.e.\ in $(0,T)\times\rr$. Since $T>0$ is arbitrary, the proof is complete.
\end{proof}

\section{The far field}\label{sect-farfield}

The goal of this section is to prove Theorems~\ref{farfield} and \ref{farfield.gamma1}: for data satisfying \eqref{initial.gamma}, far away from the origin the solution is still described by its datum, $u(t,x)\approx k_\pm|x|^{-\gamma}$, uniformly in time. This estimate is precisely what will allow us, in Sections~\ref{sect-lt1} and \ref{sect-eq1}, to pass from convergence on expanding sets to convergence in the norms to which the datum belongs.

The idea is to compare $u$ with an explicit barrier on moving sets of the form $\{x>\zeta(t)\}$. The section is organized accordingly. In Subsection~\ref{subsect.ff.comparison} we prove a comparison principle on sets of this form. In Subsection~\ref{subsect.ff.barrier} we construct the barrier. It will not be the same for $\gamma\in(0,1)$ and for $\gamma=1$, but it has the same structure in both cases, and we construct them at the same time. Finally, in Subsection~\ref{subsect.ff.proof} we first prove a far-field estimate that holds for $\gamma\in(0,1)$ and for $\gamma=1$ (Lemma~\ref{ff.general}), since the comparison argument is the same in both cases, and then we deduce the two theorems separately, by choosing the parameters of the barrier in each case.

\subsection{Comparison on a domain with a moving boundary}\label{subsect.ff.comparison}

For $\sigma\in\{1,-1\}$ and $v$ of class $C^1$ in $t$ and $C^2$ in $x$ we write
\[
\Lop_\sigma v:=v_t-v_{xx}+\sigma\,\partial_x\big(F(v)\big).
\]
For $\sigma=1$ this is the operator in \eqref{eq.main}; the case $\sigma=-1$ corresponds to the equation satisfied by $u(t,-x)$, and will be used to deal with the left tail.

We want to compare a solution $v$ with a supersolution $\bar v$ on a set $\Sigma=\{x>\zeta(t)\}$, knowing that $v\le\bar v$ on the moving boundary $x=\zeta(t)$ and, in some sense, at $t=0$. The difficulties are of two kinds. The first one concerns the equation. The difference $w=v-\bar v$ satisfies $w_t-w_{xx}+\sigma(cw)_x\le0$, with $c=(F(v)-F(\bar v))/w$ bounded; but when $1<q<2$ the function $F'$ is not Lipschitz near $0$,  and the usual proof of the maximum principle, which requires a bound on $c_x$, does not apply. The second one concerns the initial condition. The solutions to which we shall apply the result are continuous in time with values in $L^p(\rr)$, not pointwise, so at $t=0$ we can only check a condition in an integral sense.

Both difficulties are overcome by the same tool, an $L^1$ argument in the spirit of Kru\v{z}kov's proof of uniqueness \cite[Section~3]{MR267257} (see also Lemma~\ref{u:kato}). We multiply the inequality for $w$ by a smooth approximation of the Heaviside function evaluated at $w$, which only requires $c$ to be bounded, and integrate in $x$ over $(\zeta(t),\infty)$ against the weight $e^{-x}$, whose only role is to make the integral finite. The moving boundary turns out to be harmless: since $w\le0$ there, all the boundary terms vanish. What remains is a Gronwall inequality for $\int_{\zeta(t)}^\infty w_+e^{-x}\,dx$, and the initial condition we need is precisely that this integral tends to $0$ as $t\downarrow0$.

\begin{lemma}\label{ff.comparison}
    Let $\sigma\in\{1,-1\}$, let $\zeta\in C([0,\infty))\cap C^1((0,\infty))$, and let
    \[
        \Sigma:=\{(t,x):t>0,\ x>\zeta(t)\}.
    \]
    Let $v\in L^\infty(Q)$ be a classical solution of $\Lop_\sigma v=0$ in $Q$, and let $\bar v$ be of class $C^1$ in $t$ and $C^2$ in $x$ in a neighbourhood of $\overline\Sigma\cap\{t>0\}$, and bounded on $\Sigma\cap\{t<T\}$ for every $T>0$. If
    \begin{enumerate}
        \item[(i)] $\Lop_\sigma\bar v\ge0$ in $\Sigma$,
        \item[(ii)] $v(t,\zeta(t))\le\bar v(t,\zeta(t))$ for every $t>0$,
        \item[(iii)] $\displaystyle\lim_{t\downarrow0}\int_{\zeta(t)}^\infty\big(v(t,x)-\bar v(t,x)\big)_+e^{-x}\,dx=0$,
    \end{enumerate}
    then $v\le\bar v$ in $\Sigma$.
\end{lemma}

\begin{proof}
Let $w:=v-\bar v$ and $h:=F(v)-F(\bar v)$, so that, by (i), $w_t-w_{xx}+\sigma h_x\le0$ in $\Sigma$. We have to prove that $w_+=0$.

Fix $T>0$. Since $v$ and $\bar v$ are bounded on $\Sigma\cap\{t<T\}$, there is a constant $L$, the Lipschitz constant of $F$ on the range of $v$ and $\bar v$, such that $|h|\le L|w|$ there. Let $S\in C^\infty(\rr)$ be nondecreasing, with $S=0$ on $(-\infty,0]$ and $S=1$ on $[1,\infty)$, and for $\eta>0$ let $S_\eta(r):=S(r/\eta)$ and $\Phi_\eta(r):=\int_0^rS_\eta$. Thus $0\le S_\eta\le1$ and $S_\eta'\ge0$, and $\Phi_\eta$ is a smooth convex approximation of $r_+$, with $0\le\Phi_\eta(r)\le r_+$ and $\Phi_\eta(r)\to r_+$ as $\eta\downarrow0$. We multiply the inequality for $w$ by $S_\eta(w)\ge0$ and use
\begin{gather*}
    S_\eta(w)w_t=\partial_t\Phi_\eta(w),\quad S_\eta(w)w_{xx}=\partial_{xx}\Phi_\eta(w)-S_\eta'(w)w_x^2,\\
    S_\eta(w)h_x=\partial_x\big(S_\eta(w)h\big)-S_\eta'(w)w_xh .
\end{gather*}
We obtain, in $\Sigma\cap\{t<T\}$,
\[
    \partial_t\Phi_\eta(w)-\partial_{xx}\Phi_\eta(w)+\sigma\,\partial_x\big(S_\eta(w)h\big)\le S_\eta'(w)\big(\sigma\,w_xh-w_x^2\big) .
\]
Since $S_\eta'\ge0$ and $|w_xh|\le L|w||w_x|\le w_x^2+\frac{L^2}4w^2$, the right-hand side is at most $\frac{L^2}4S_\eta'(w)w^2$. Moreover, since $S_\eta'(w)=\frac1\eta S'(w/\eta)$ vanishes unless $0<w<\eta$,
\[
    S_\eta'(w)\,w^2\le\frac1\eta\,\|S'\|_{L^\infty(\rr)}\,\eta^2=\|S'\|_{L^\infty(\rr)}\,\eta,
\]
and therefore
\[
    \partial_t\Phi_\eta(w)-\partial_{xx}\Phi_\eta(w)+\sigma\,\partial_x\big(S_\eta(w)h\big)\le\frac{L^2}4\,S_\eta'(w)\,w^2\le c_\eta:=\frac{L^2}4\,\|S'\|_{L^\infty(\rr)}\,\eta .
\]

Next we integrate in $x$ against the weight $\psi_n(x):=e^{-x}\chi(x/n)$, where $\chi\in C^\infty_c(\rr)$ is supported in $[-2,2]$, $0\le\chi\le1$ and $\chi=1$ on $[-1,1]$; thus $|\psi_n'|+|\psi_n''|\le c_ne^{-x}$, with $c_n\to2$ as $n\to\infty$. Let
\[
    Y(t):=\int_{\zeta(t)}^\infty\Phi_\eta(w(t,x))\,\psi_n(x)\,dx,\qquad 0<t<T .
\]
Since $\psi_n$ vanishes for $x\ge2n$, $Y(t)$ is an integral over a bounded interval, whose left endpoint $\zeta(t)$ is of class $C^1$, of a function of class $C^1$ in $t$; by Leibniz's rule,
\[
    Y'(t)=\int_{\zeta(t)}^\infty\partial_t\Phi_\eta(w)\,\psi_n\,dx-\zeta'(t)\,\Phi_\eta\big(w(t,\zeta(t))\big)\,\psi_n(\zeta(t)) .
\]
By (ii), $w\le0$ on $x=\zeta(t)$, so that $\Phi_\eta(w)$, $\partial_x\Phi_\eta(w)=S_\eta(w)w_x$ and $S_\eta(w)h$ all vanish on $x=\zeta(t)$. In particular the last term in $Y'(t)$ vanishes. For the first one we use the inequality for $\Phi_\eta(w)$, multiplying it by $\psi_n\ge0$, integrating over $(\zeta(t),\infty)$, and integrating by parts in $x$. No boundary terms appear, since they vanish at infinity because $\psi_n$ has compact support, and at $x=\zeta(t)$, as we have just checked. This gives, for $0<t<T$,
\begin{align*}
    Y'(t)&\le\int_{\zeta(t)}^\infty\big(\Phi_\eta(w)\,\psi_n''+\sigma\,S_\eta(w)h\,\psi_n'+c_\eta\psi_n\big)\,dx\\
    &\le\int_{\zeta(t)}^\infty\big(\Phi_\eta(w)|\psi_n''|+L\,w_+|\psi_n'|+c_\eta\psi_n\big)\,dx,
\end{align*}
where in the last step we have used that $S_\eta(w)|h|\le L\,w_+$, because $0\le S_\eta\le1$, $S_\eta(w)=0$ if $w\le0$, and $|h|\le L|w|$.
We integrate on $(t_1,t_2)$ and let $\eta\downarrow0$ and then $n\to\infty$. Setting $Z(t):=\int_{\zeta(t)}^\infty w_+(t,x)\,e^{-x}dx$, we obtain
\[
    Z(t_2)\le Z(t_1)+2(1+L)\int_{t_1}^{t_2}Z(t)\,dt\qquad\text{for }0<t_1<t_2<T,
\]
and Gronwall's lemma gives $Z(t_2)\le e^{2(1+L)(t_2-t_1)}Z(t_1)$. Letting $t_1\downarrow0$ and using (iii), $Z\equiv0$ on $(0,T)$. Since $w$ is continuous and $T$ is arbitrary, $w\le0$ in $\Sigma$.
\end{proof}

\subsection{The barrier}\label{subsect.ff.barrier}

Throughout this subsection $q>1$, $\gamma\in(0,1]$ and $\delta:=\gamma(q-1)$.

On the inner boundary $x=\zeta(t)$ we know nothing about $u$ beyond the decay estimates~\eqref{thm1.decay} and \eqref{thm2.decay}, so the barrier has to be at least as large as the decay rate there. Inside the set it has to be a supersolution, and at the initial time it has to lie above the datum. To make a good guess for it, take $u_0=kx^{-\gamma}$ and expand the solution to first order in time: $u(t,x)\approx u_0(x)+t\,u_t(0,x)$, with $u_t(0,x)=u_0''(x)-\big(F(u_0)\big)'(x)$. Since $\gamma q+1=\gamma+1+\delta$, this gives
\begin{equation}\label{ff.formal}
    x^\gamma u(t,x)-k\approx k\big[\gamma(\gamma+1)\,t\,x^{-2}+q\gamma|k|^{q-1}\,t\,x^{-1-\delta}\big].
\end{equation}
The first term is the correction due to diffusion, and the second one the correction due to convection. Notice that the relation between space and time in each of them is the one dictated by the homogeneity of the corresponding problem, which confirms that we are on the right track. For $k\ne0$, both corrections are bounded on $\{x\ge Lt^{\theta}\}$, uniformly in $t\ge1$.

The barrier is built from this expansion. It has the form
\begin{equation}\label{ff.bar}
    \bar v(t,x)=(k+\eps)x^{-\gamma}+E\,t\,x^{-\gamma-2}+D(t)\,t\,x^{-\gamma q-1}+B(t)\,\zeta(t)\,x^{-\gamma-1},\qquad x>\zeta(t),
\end{equation}
where $x=\zeta(t)$ is the inner boundary and $E$, $D(t)$, $B(t)$ are positive. The first term is the datum, shifted by $\eps$ so that it lies above $u_0$. The next two are the corrections in \eqref{ff.formal}, with coefficients large enough to make $\bar v$ a supersolution, since they have to absorb what the operator produces when acting on the datum and on the other terms. The last term, the \emph{boundary term}, is introduced so that the barrier lies above the solution on the moving boundary; it has to be unimportant farther inside. We can write it as $\big(B(t)\zeta(t)/x\big)\,x^{-\gamma}$, that is, $B(t)\zeta(t)/x$ times the datum profile $x^{-\gamma}$. This factor equals $B(t)$ on the inner boundary, where $B(t)$ is chosen so that the barrier dominates the decay rate, and the boundary term is negligible with respect to the datum when $|x|\gg B(t)\zeta(t)$.

The barrier will be used both for $\gamma\in(0,1)$ and for $\gamma=1$, and the two cases differ in the inner boundary and in the decay rate that the barrier has to dominate on it. For $\gamma\in(0,1)$ the inner boundary is at distance $(1+t)^\theta$ and the decay rate is $C(1+t)^{-\gamma\theta}$, where $C$ is the constant in \eqref{thm1.decay}. The coefficient $B$ is then a constant, and the boundary term, like the two corrections, is small for $|x|\gg t^\theta$, which is the scale of Theorem~\ref{farfield}. For $\gamma=1$ the inner boundary is at the scale $m(t)$ of \eqref{m.def}, which is $(1+t)^\theta$ if $q>2$ and carries a logarithm if $q\le2$. In terms of this scale \eqref{thm2.decay} reads $\|u(t)\|_{L^\infty(\rr)}\le C\log(2+t)/m(t)$ for every $q>1$, while the datum $|x|^{-1}$ has size $1/m(t)$ there. Hence we need $B(t)\approx\log(2+t)$ instead of a constant. Then the boundary term is negligible only for $|x|\gg B(t)\zeta(t)\approx m(t)\log(2+t)$: the far field begins one logarithm beyond the scale $m(t)$ of the inner region, as stated in Theorem~\ref{farfield.gamma1}.

To treat both cases at once we take
$  \zeta(t)= {\Lambda m(t)} $ and  $Cb(t)m^{-\gamma}(t)$  as an upper bound for the decay rate,
where $\Lambda$ is a large constant. Concretely, we take $(m,b)=\big((1+t)^\theta,1\big)$ if $\gamma\in(0,1)$, and, if $\gamma=1$, $m$ as in \eqref{m.def} and $b(t)=\log(2+t)$. It is easy to check that, in both cases,
\begin{equation}\label{m.props}
    \begin{gathered}
        m,b\in C^1([0,\infty))\ \text{nondecreasing},\quad m\geq m_*>0,\qquad b\ge b_*>0,\\
        \mu_1:=\sup_{t\ge0}\frac t{m^2(t)}<\infty,\quad
        \mu_2:=\sup_{t\ge0}\frac{t\,b^{q-1}(t)}{m^{1+\delta}(t)}<\infty,\\
        \mu_3:=\inf_{t\ge0}m'(t)\,m(t)>0,\quad
        \mu_4:=\inf_{t\ge0}\frac{m'(t)\,m^\delta(t)}{b^{q-1}(t)}>0 .
    \end{gathered}
\end{equation}
These are the only properties of $m$ and $b$ that we need to prove the barrier lemma below. The two suprema say that the corrections in \eqref{ff.formal} stay bounded on the inner boundary (the factor $b^{q-1}(t)$ accounts for the size $b(t)$ of the barrier there), and the two infima that the boundary moves fast enough for the boundary term, which grows in time through $\zeta(t)$, to be a supersolution.

\begin{lemma}\label{ff.barrier}
    Let $q>1$, $\gamma\in(0,1]$,  $k\in\rr$, $C>0$, and let $m$, $b$ satisfy \eqref{m.props}. There are $E=E(|k|,\gamma)$ and  $C_0>0$ such that for every  $\Lambda $ large enough there are positive nondecreasing functions $B$, $D$ on $[0,\infty)$, satisfying
        \begin{equation}
          \label{eq:BandD}
        B(t)\le C_0\Lambda^\gamma b(t), \quad D(t)\le C_0 \Lambda^\delta b^{q-1}(t) \quad \text{for every}\  t\ge 0,   
        \end{equation}
      such that for every $\eps\in(0,1]$ the function $\bar v$ in \eqref{ff.bar}, with $\zeta(t)=\Lambda m(t)$, has the following properties, where $\Sigma:=\{(t,x):t>0,\ x>\zeta(t)\}$:
    \begin{enumerate}
        \item[(a)] $\bar v$ lies above the datum: $\bar v\ge(k+\eps)x^{-\gamma}$ in $\Sigma$;
        \item[(b)] $\bar v$ is of class $C^1$ in $t$ and $C^2$ in $x$ in a neighbourhood of $\overline\Sigma\cap\{t>0\}$, and bounded on $\Sigma\cap\{t<T\}$ for every $T>0$;
        \item[(c)] $\bar v$ dominates the decay rate on the boundary: $\bar v(t,\zeta(t))\ge C\,b(t)\,m^{-\gamma}(t)$ for every $t>0$;
        \item[(d)] $\bar v$ is a supersolution: $\Lop_\sigma\bar v\ge0$ in $\Sigma$ for $\sigma=\pm1$;
    \end{enumerate}
\end{lemma}

\begin{proof}
We look for $\bar v$ of the form \eqref{ff.bar}, with $\zeta(t)=\Lambda m(t)$, $E>0$ and nonnegative functions $B$, $D$ of class $C^1$ still to be chosen for $\Lambda$ large enough. For any such choice, (a) is clear, since all the terms of $\bar v-(k+\eps)x^{-\gamma}$ are nonnegative, and (b) holds because $\bar v$ is of class $C^1$ in $t$ and $C^2$ in $x$ for $x>0$, its coefficients are bounded for $t<T$, and $x>\Lambda m_*$ in $\Sigma$. Then the proof goes as follows. Property (c) determines $B$. Property (d) determines $E$ and $D$.   The functions $B$ and $D$ that we obtain are of class $C^1$ because $b$ and $m$ are.

 At $x=\zeta(t)$ the second and third terms in \eqref{ff.bar} are nonnegative and the last one equals $B(t)\zeta^{-\gamma}(t)$, so that $\bar v(t,\zeta(t))\ge\big(B(t)-|k|\big)\zeta^{-\gamma}(t)$. Since $\zeta(t)=\Lambda m(t)$, the choice
\[
    B(t):= |k|+C\Lambda^\gamma b(t), \ 
\]
gives \[\bar v(t,\zeta(t))\ge  C\Lambda^\gamma b(t)\,\zeta^{-\gamma}(t)=Cb(t)m^{-\gamma}(t),\] which is (c). Note that $B $ satisfies \eqref{eq:BandD} for $\Lambda$ large enough and  $B$ is nondecreasing, because $b$ is. 

 Write $\bar v=P_0+P_1+P_2+P_3$ in the order of \eqref{ff.bar}. We shall choose $\Lambda$ so large that the two corrections satisfy
\begin{equation}\label{ff.corr}
    E\,t\,x^{-2}\le |k|+1,\quad D(t)\,t\,x^{-1-\delta}\le |k|+1\quad\text{in }\Sigma.
\end{equation}
 Then, since $|k+\eps|\le |k|+1$ and $\zeta\le x$ in $\Sigma$,
\[
    |\bar v(t,x)|\le K(t)\,x^{-\gamma}\quad\text{in }\Sigma,\qquad K(t):=3(|k|+1)+B(t),
\]
so that $|F'(\bar v)|=q|\bar v|^{q-1}\le a(t)\,x^{-\delta}$ with $a(t):=qK^{q-1}(t)$, and
\[
    \Lop_\sigma\bar v=\sum_{i=0}^3\big(\partial_tP_i-\partial_{xx}P_i\big)+\sigma F'(\bar v)\,\partial_x\bar v\ge\sum_{i=0}^3J_i,\quad J_i:=\partial_tP_i-\partial_{xx}P_i-a(t)\,x^{-\delta}\,|\partial_xP_i| .
\]
For $i=1,2,3$, $P_i$ is a nonnegative multiple of $x^{-\lambda_i}$, with $\lambda_1=\gamma+2$, $\lambda_2=\gamma q+1=\gamma+1+\delta$ and $\lambda_3=\gamma+1$, so that $\partial_{xx}P_i=\lambda_i(\lambda_i+1)P_i/x^2$ and $|\partial_xP_i|=\lambda_iP_i/x$; moreover $\partial_tP_1=Ex^{-\gamma-2}$ and, recalling that $B$ is nondecreasing, $\partial_tP_3\ge P_3\zeta'/\zeta$; likewise $\partial_tP_2\ge D(t)x^{-\gamma q-1}$ if $D$ is nondecreasing, as will be the case. Hence, using $|k+\eps|\le |k|+1$ in $J_0$ and $x\ge\zeta$ in $J_3$,
\[
    \begin{aligned}
        J_0&\ge-\gamma(\gamma+1)(|k|+1)\,x^{-\gamma-2}-\gamma(|k|+1)\,a(t)\,x^{-\gamma q-1},\\
        J_1&=Ex^{-\gamma-2}\big[1-(\gamma+2)(\gamma+3)\,t\,x^{-2}-(\gamma+2)\,a(t)\,t\,x^{-1-\delta}\big],\\
        J_2&\ge D(t)\,x^{-\gamma q-1}\big[1-(\gamma q+1)(\gamma q+2)\,t\,x^{-2}-(\gamma q+1)\,a(t)\,t\,x^{-1-\delta}\big],\\
        J_3&\ge\frac{P_3}{\zeta}\Big[\zeta'-\frac{(\gamma+1)(\gamma+2)}{\zeta}-\frac{(\gamma+1)\,a(t)}{\zeta^\delta}\Big] .
    \end{aligned}
\]
There is a $\kappa_0$ independent of $\Lambda$ such that
\begin{equation}\label{ff.Ka}
      a(t)\le \kappa_0\Lambda^\delta b^{q-1}(t)
\end{equation}
for  $\Lambda$ large enough.
Hence, by \eqref{m.props} and since $x>\Lambda m(t)$ in $\Sigma$, we have
\begin{equation}\label{ff.tx}
    t\,x^{-2}\le\frac {\mu_1}{\Lambda^2},\quad a(t)\,t\,x^{-1-\delta}\le\frac{\kappa_0\mu_2}{\Lambda}\quad\text{in }\Sigma.
\end{equation}
Then for $\Lambda$  large enough 
\[J_1+J_2\ge\frac12E\,x^{-\gamma-2}+\frac12D(t)\,x^{-\gamma q-1},\] which cancels the two negative terms in $J_0$ if
\[
    E:=2\gamma(\gamma+1)(|k|+1),\qquad D(t):=2\gamma(|k|+1)\,a(t).
\]
Notice that with this choice of $D$ and $E$ \eqref{ff.corr} holds possibly taking  a larger $\Lambda$.
Moreover, $D$ is nondecreasing, since $K$ is, and satisfies \eqref{eq:BandD}. 

Finally, the bracket in $J_3$ is nonnegative if
\begin{equation}\label{ff.ell}
    \zeta'(t)\,\zeta(t)\ge2(\gamma+1)(\gamma+2),\quad\zeta'(t)\,\zeta^\delta(t)\ge2(\gamma+1)\,a(t)\quad\text{for }t>0 .
\end{equation}

  It remains to satisfy   \eqref{ff.ell}. Since $\zeta=\Lambda m$ and $\zeta'=\Lambda m'$, \eqref{m.props} and \eqref{ff.Ka} give, for $t>0$,
\[
    \zeta'(t)\,\zeta(t)\ge\Lambda^2\mu_3,\qquad
    \zeta'(t)\,\zeta^\delta(t)\ge\Lambda^{1+\delta}\mu_4\,b^{q-1}(t)\ge\frac{\Lambda\mu_4}{\kappa_0}\,a(t) .
\]
Hence  \eqref{ff.ell} holds taking $\Lambda$ large enough. This concludes the proof of (d).
\end{proof}

\subsection{Proof of Theorems~\ref{farfield} and \ref{farfield.gamma1}}\label{subsect.ff.proof}

The comparison argument is the same for $\gamma\in(0,1)$ and for $\gamma=1$, and it only involves the functions $m$ and $b$ of the barrier. We first carry it out for general $m$ and $b$, which gives a far-field estimate (Lemma~\ref{ff.general}), and then we prove the two theorems separately, choosing $m$ and $b$ in each case. In what follows $\delta:=\gamma(q-1)$.

\begin{lemma}\label{ff.general}
    Let $\gamma\in(0,1]$, $q>1$, $u_0$ satisfying \eqref{initial.gamma}, and $m$ and $b$ satisfying \eqref{m.props} and
    \begin{equation}\label{ff.decay}
        \|u(t)\|_{L^\infty(\rr)}\le C\,b(t)\,m^{-\gamma}(t)\quad\text{for all }t\ge0,
    \end{equation}
    for some $C>0$. Then, for every $\eps\in(0,1]$ there is $K_\eps\ge1$ such that
    \begin{equation}\label{ff.quant}
        \big|\,|x|^\gamma u(t,x)-k_\pm\big|\le\eps+K_\eps\Big(\frac t{x^2}+\frac{b^{q-1}(t)\,t}{|x|^{1+\delta}}+\frac{b(t)\,m(t)}{|x|}\Big)\quad\text{for }t>0,\ \pm x\ge K_\eps\,m(t).
    \end{equation}
\end{lemma}

\begin{proof}
 We prove the upper estimate for $u$.   The other estimates are obtained analogously using    $-u$, $\tilde u$ and $-\tilde u$, where $\tilde u(t,x):=u(t,-x)$.
  
    Fix $\eps\in(0,1]$ and let $R_\eps\ge1$ be such that $\big|\,|x|^\gamma u_0(x)-k_\pm\big|\le\eps$ for a.e.\ $\pm x\ge R_\eps$. We take $\Lambda_\eps\ge R_\eps/m_*$ so large that Lemma~\ref{ff.barrier} applies with $k=k_+$. 
    
Let $\bar v$ be the barrier \eqref{ff.bar} with $k=k_+$, and let us check the assumptions of Lemma~\ref{ff.comparison} for $v=u$ and $\sigma=1$. The regularity of $\bar v$ and (i) are Lemma~\ref{ff.barrier}(b) and (d). Condition (ii) holds because $u(t,\zeta(t))\le C\,b(t)\,m^{-\gamma}(t)\le\bar v(t,\zeta(t))$, by \eqref{ff.decay} and Lemma~\ref{ff.barrier}(c). As for~(iii), the datum $u_0$ satisfies $u_0\le(k_++\eps)x^{-\gamma}\le\bar v(t,\cdot)$ a.e.\ on $\{x>\zeta(t)\}$,  using $\zeta(t)=\Lambda_\eps m(t)\ge R_\eps$ and   
Lemma~\ref{ff.barrier}(a). Hence
    \[
        \big(u(t)-\bar v(t)\big)_+\le|u(t)-u_0|
    \]
    there, and, since $u\in C([0,\infty);L^{2/\gamma}(\rr))$, the integral in  (iii) tends to $0$ by H\"older's inequality. Lemma~\ref{ff.comparison} gives $u\le\bar v$ in $\Sigma$, which is     \[
        x^\gamma u(t,x)-k_+\le\eps+E\,t\,x^{-2}+D(t)\,t\,x^{-1-\delta}+B(t)\,\frac{\zeta(t)}x\qquad\text{for }t>0,\ x>\zeta(t),
    \]
    from which  
 \eqref{ff.quant} follows using \eqref{eq:BandD} after taking $K_\eps>\Lambda_\eps$ large enough. 
\end{proof}

We now deduce Theorem~\ref{farfield} from Lemma~\ref{ff.general}. In this case the coefficient of the boundary term of the barrier can be taken constant, $b\equiv1$, and the far field starts at the scale $t^\theta$ of the inner region.

\begin{proof}[Proof of Theorem~\ref{farfield}]
    Let $\gamma\in(0,1)$. We apply Lemma~\ref{ff.general} with $m(t)=(1+t)^\theta$ and $b\equiv1$, so that \eqref{ff.decay} is \eqref{thm1.decay}.

    Let us prove \eqref{thm.ff}. Fix $\eps\in(0,1]$ and let $K_\eps$ be given by Lemma~\ref{ff.general}. Let 
     $|x|\ge L(1+t^\theta)\geq L\,m(t)$. If $L\ge K_\eps$  then \eqref{ff.quant} applies. Moreover, the definitions of $\mu_1$ and $\mu_2$ in \eqref{m.props}  give that $t/x^2\le\mu_1L^{-2}$ and $t/|x|^{1+\delta}\le\mu_2L^{-1-\delta}$. Hence
    \[
        \sup_{t>0}\ \sup_{\pm x\ge L(1+t^\theta)}\big|\,|x|^\gamma u(t,x)-k_\pm\big|\le\eps+K_\eps\big(\mu_1L^{-2}+\mu_2L^{-1-\delta}+L^{-1}\big).
    \]
  Letting first $L\rightarrow\infty$ and then $\eps\rightarrow 0$  proves \eqref{thm.ff}.

    To conclude, we have to prove~\eqref{thm.env}. By \eqref{thm.ff} there is $L$ such that
    \[
        |x|^\gamma|u|\le   \max\{|k_+|,|k_-|\}+1\quad\textup{on }\{|x|\ge L(1+t^\theta)\}.
    \]
    There $1+t^\theta+|x|\le2|x|$, so that $\big(1+t^\theta+|x|\big)^\gamma|u|\le(2|x|)^\gamma|u|\le \mathcal{N}$. On the complement, $1+t^\theta+|x|\le2(1+L)(1+t)^\theta$, so that
    \[
        \big(1+t^\theta+|x|\big)^\gamma|u|\le \mathcal{N}(1+t)^{\gamma\theta}\|u(t)\|_{L^\infty(\rr)}\le \mathcal{N}
    \]
    by \eqref{thm1.decay}.
\end{proof}

The proof of Theorem~\ref{farfield.gamma1} follows the same lines, with a choice of $m$ and $b$ that accounts for the logarithm: the coefficient of the boundary term grows like $\log(2+t)$, and this moves the beginning of the far field to the scale $m(t)\log(2+t)$.

\begin{proof}[Proof of Theorem~\ref{farfield.gamma1}]As in Lemma \ref{ff.general} we prove the estimate for $u$ and $x>0$. 
    In this case  $\delta=q-1$. We apply Lemma~\ref{ff.general} with $m$ as in \eqref{m.def} and $b(t):=\log(2+t)$. In terms of $m$, \eqref{thm2.decay} reads $\|u(t)\|_{L^\infty(\rr)}\le C\log(2+t)/m(t)$ for every $q>1$, so that \eqref{ff.decay} holds.

    Let us prove \eqref{ff.gamma1}. Fix $\eps\in(0,1]$ and let $K_\eps$ be given by Lemma~\ref{ff.general}. Let $L\ge K_\eps/\log(2)$ and $x\ge L\,m(t)\log(2+t)$. Then $x\ge K_\eps m(t)$, so that \eqref{ff.quant} applies. 
    The same argument as in the proof of  Theorem~\ref{farfield}
   shows \eqref{ff.gamma1}. 

    We now proceed with the proof of \eqref{env.gamma1}. Taking $K$ such that  $|x|\ge K m(t)$, \eqref{ff.quant} gives
    \[
        |x|\,|u|\le \max\{|k_+|,|k_-|\}+1+K\Big(\frac t{x^2}+\frac{b^{q-1}(t)\,t}{|x|^q}+\frac{b(t)\,m(t)}{|x|}\Big)\le \mathcal{M}+\mathcal{M}\,\frac{m(t)\log(2+t)}{|x|},
    \]
    since, by \eqref{m.props}, the first two terms in the parenthesis are bounded for $|x|\ge Km(t)$.  In the region $0<|x|<Km(t)$ we use \eqref{thm2.decay}: \[
    |u(t,x)|\le C\frac{\log(2+t)}{m(t)}\le CK^2\,\frac{m(t) \log(2+t)}{x^2}.\] This proves~\eqref{env.gamma1}.
\end{proof}

\section{The case \texorpdfstring{$\gamma\in(0,1)$}{gamma in (0,1)}: proof of Theorem~\ref{first.term.gamma<1}}\label{sect-lt1}

Throughout this section $\gamma\in(0,1)$, $u_0$ satisfies \eqref{initial.gamma}, and
\[\kappa:=\esssup_{x\in\rr}(1+|x|)^\gamma|u_0(x)|.\]
We first prove   the convergence in the inner region
\begin{equation}\label{loc.conv}
\lim_{t\to\infty}\int_{|\xi|\le R}\big|t^{\gamma\theta}u(t,t^\theta\xi)-U_\gamma(\xi)\big|\,d\xi=0\qquad\text{for every }R>0 ,
\end{equation}
and then, in Subsection~\ref{subsect.upgrade}, we upgrade \eqref{loc.conv} to \eqref{thm1.conv}--\eqref{thm1.conv.p} by means of Theorem~\ref{farfield}.

We  record, in the following lemma, whose proof is standard, the behaviour of the heat equation with data \eqref{initial.gamma}; related results can be found in \cite{KaUg87,Herraiz-1999}.

\begin{lemma}\label{decay.linear.heat}
Let $\gamma\in(0,1)$ and let $u_0\in L^\infty(\rr)$ satisfy \eqref{initial.gamma}. Then
\[\lim _{t\rightarrow\infty}t^{\gamma/2}(G_t\ast u_0)(t^{1/2}\xi) = (G_1\ast U_0)(\xi)\]
 uniformly for $\xi$ in compact sets.
\end{lemma}

\subsection{The case \texorpdfstring{$q>1+\frac1\gamma$}{q>1+1/gamma}}\label{subsect.supercritical}

Here $\theta=\frac12$, so that $\|u(s)\|_{L^\infty(\rr)}\le C(1+s)^{-\gamma/2}$ by \eqref{thm1.decay}. Let $a:=\frac{q\gamma-1}2$. By \eqref{interp.lorentz} and \eqref{est.p.infty},
\[
    \|F(u(s))\|_{L^1(\rr)}=\|u(s)\|_{L^q(\rr)}^q\le C\|u(s)\|_{L^{1/\gamma,\infty}(\rr)}^{1/\gamma}\|u(s)\|_{L^\infty(\rr)}^{q-\frac1\gamma}\le C(1+s)^{-a},
\]
and $\|F(u(s))\|_{L^\infty(\rr)}\le C(1+s)^{-q\gamma/2}$. By \eqref{duhamel} and the bounds $\|\partial_xG_\tau\|_{L^\infty(\rr)}\le C\tau^{-1}$ and $\|\partial_xG_\tau\|_{L^1(\rr)}\le C\tau^{-1/2}$, used for $s<\frac t2$ and $s>\frac t2$, respectively,
\begin{align*}
    t^{\gamma/2}\|u(t)-G_t\ast u_0\|_{L^\infty(\rr)}&\le Ct^{\frac\gamma2-1}\int_0^{t/2}(1+s)^{-a}\,ds+Ct^{\frac\gamma2}\int_{t/2}^t(t-s)^{-\frac12}(1+s)^{-\frac{q\gamma}2}\,ds \\
    &\lesssim t^{\frac\gamma2-1}\max\{1,\log(1+t),t^{1-a}\}+t^{(1+\gamma-q\gamma)/2}.
\end{align*}
Both terms tend to $0$ as $t\rightarrow \infty$, since $\gamma<2$ and $a>\gamma/2$. Lemma~\ref{decay.linear.heat}
shows $t^{\gamma/2}u(t,t^{1/2}\xi)\to U_\gamma(\xi)$ uniformly on compact sets and
 \eqref{loc.conv}.

\subsection{The case \texorpdfstring{$1<q\le1+\frac1\gamma$}{1<q<=1+1/gamma}}\label{subsect.critical}
Now $\theta=\beta$ and $\gamma\theta=\alpha$. For $\lambda>1$ let
\[
    u_\lambda(t,x):=\lambda^\gamma u(\lambda^{1/\beta}t,\lambda x),\qquad u_{0,\lambda}(x):=\lambda^\gamma u_0(\lambda x).
\]
Since $\gamma q+1=\gamma+\frac1\beta$, $u_\lambda$ solves
\begin{equation}\label{visc.lambda}
    \partial_tu_\lambda+\partial_x(F(u_\lambda))=\eps_\lambda\partial_{xx}u_\lambda,\qquad\eps_\lambda:=\lambda^{(q-1)\gamma-1}\in(0,1],
\end{equation}
with initial datum $u_{0,\lambda}$. If $q=1+\frac1\gamma$, then $\beta=\frac12$, $\alpha=\frac\gamma2$ and $\eps_\lambda=1$, so that $u_\lambda$ is the solution of \eqref{eq.main} with initial datum $u_{0,\lambda}$; if $q<1+\frac1\gamma$, then $\eps_\lambda\to0$ as $\lambda\to\infty$. Moreover $|u_{0,\lambda}(x)|\le\kappa|x|^{-\gamma}$, and $u_{0,\lambda}\to U_0$ in $L^1_{\textup{loc}}(\rr)$ as $\lambda\to\infty$.

\begin{lemma}\label{uniform.gamma.subcritical}
    There are constants $C$ and $C_p$, independent of $\lambda>1$ and $t>0$, such that
    \begin{enumerate}
        \item[(a)] $\|u_\lambda(t)\|_{L^\infty(\rr)}\le Ct^{-\alpha}$;
        \item[(b)] $\int_E|u_\lambda(t,x)|\,dx\le\frac{2^\gamma\kappa}{1-\gamma}|E|^{1-\gamma}$ for every measurable set $E\subset\rr$;
        \item[(c)] $\int_0^T\!\!\int_E|u_\lambda|^p\,dx\,dt\le C_p\,|E|^{1-\gamma}\,T^{1-\alpha(p-1)}$ for every measurable $E\subset\rr$, $T>0$ and $p\in[1,q+\frac1\gamma)$.
    \end{enumerate}
\end{lemma}

\begin{proof}
Parts (a) and (b) follow from the second bound in \eqref{decay} and from \eqref{local.L1.est}, written for $u$, through the change of variables defining $u_\lambda$, since $\alpha=\gamma\beta$. Part (c) follows from $\int_E|u_\lambda(t)|^p\le\|u_\lambda(t)\|_{L^\infty(\rr)}^{p-1}\int_E|u_\lambda(t)|$, (a), (b), and $\alpha(p-1)<1$ for $p<1+\frac1\alpha=q+\frac1\gamma$.
\end{proof}

Let $0<\tau<T$, $R>0$, and $\psi\in C^\infty_c(\rr)$ with $0\le\psi\le1$, $\psi=1$ on $(-R,R)$. Multiplying \eqref{visc.lambda} by $u_\lambda\psi$
  we get
\[
\eps_\lambda\int_\tau^T\!\!\int_\rr|\partial_xu_\lambda|^2\psi\,dx\,dt\le\frac12\int_\rr u_\lambda(\tau)^2\psi\,dx+\int_\tau^T\!\!\int_\rr\Big(\frac{\eps_\lambda}2u_\lambda^2|\psi''|+|u_\lambda|^{q+1}|\psi'|\Big)dx\,dt,
\]
which is bounded independently of $\lambda$ by Lemma~\ref{uniform.gamma.subcritical}(a),(b). Hence $\eps_\lambda^{1/2}\partial_xu_\lambda$ is bounded in $L^2((\tau,T)\times(-R,R))$, and $u_\lambda$ is bounded in $L^\infty((\tau,T)\times\rr)$.

Classical arguments as in \cite{endalJMPA}, \cite{CazacuIgnatPazoto} show that
every sequence $\lambda_j\to\infty$ has a subsequence, not relabelled, such that $u_{\lambda_j}\to U$ in $L^1_{\textup{loc}}(Q)$ and a.e.\ in $Q$. Passing to the limit in \eqref{visc.lambda}, $U$ solves~\eqref{eq.main} in $\mathcal D'(Q)$ if $q=1+\frac1\gamma$; passing to the limit in the entropy inequalities,  $U$ satisfies \eqref{ineg.k} if $q<1+\frac1\gamma$. Moreover, by Lemma~\ref{uniform.gamma.subcritical}(a) and this convergence, hypotheses \textup{(H1)} and \textup{(H2)} of Proposition~\ref{n:fixedtime} hold on every slab $(\tau,T)\times\rr$. Hence $U$ has a representative in $C((0,\infty);L^1_{\textup{loc}}(\rr))$,
\begin{equation}\label{n:everyt-sub}
u_{\lambda_j}(t)\to U(t)\quad\text{in }L^1_{\textup{loc}}(\rr)\qquad\text{for every }t>0,
\end{equation}
and $\|U(t)\|_{L^\infty(\rr)}\le Ct^{-\alpha}$ for every $t>0$, where $\alpha<1/q$.

\begin{proposition}\label{u:trace.sub}
    The limit $U$ satisfies \eqref{initial.entropy} with initial datum $U_0$.
\end{proposition}
\begin{proof}
Fix $R>0$, set $R'=R+1$, and let $\varrho>0$. Choose $\varsigma\in C^2_c(\rr)$ such that $\|U_0-\varsigma\|_{L^1(-R',R')}<\varrho$, and let $\phi\in C^2_c((-R',R'))$ with $0\le\phi\le1$ and $\phi=1$ on $[-R,R]$. The function $\varsigma$ solves \eqref{visc.lambda} up to the error $(F(\varsigma))'-\eps_\lambda\varsigma''$, and Remark~\ref{kato.function} applies to $u_\lambda$ and $\varsigma$ with $\nu=\eps_\lambda$.
 Letting $t_1\downarrow0$ and using Lemma~\ref{uniform.gamma.subcritical}(c), we get
\begin{multline*}
\int_\rr|u_\lambda(t)-\varsigma|\phi\,dx\le\int_\rr|u_{0,\lambda}-\varsigma|\phi\,dx\\
+\int_0^t\!\!\int_\rr\big(|F(u_\lambda)-F(\varsigma)|\,|\phi'|+\eps_\lambda|u_\lambda-\varsigma|\,|\phi''|\big)\,dx\,ds+C_\varsigma\,t,
\end{multline*}
where
$ C_\varsigma:=q\|\varsigma\|_{L^\infty(\rr)}^{q-1}\|\varsigma'\|_{L^1(\rr)}+\|\varsigma''\|_{L^1(\rr)}.$
By Lemma~\ref{uniform.gamma.subcritical}(c) with $p=q$ and $p=1$, applied on $E=(-R',R')$, the double integral is bounded by $\mathcal G\big(t^{1-\alpha(q-1)}+t\big)$, where $\mathcal G$ depends on $\varsigma$, $R'$, $\gamma$, $q$ and $\kappa$, but not on $\lambda$ or $t$. We let $\lambda=\lambda_j\to\infty$, using \eqref{n:everyt-sub}, $u_{0,\lambda}\to U_0$ in $L^1(-R',R')$ and $\|U_0-\varsigma\|_{L^1(-R',R')}<\varrho$, and conclude that
\[
\int_{-R}^R|U(t)-U_0|\,dx\le2\varrho+(\mathcal G+C_\varsigma)\big(t^{1-\alpha(q-1)}+t\big)\qquad\text{for all }t>0 .
\]
Since $\alpha(q-1)<1$, the $\limsup$ as $t\downarrow0$ of the left-hand side is at most $2\varrho$, and $\varrho>0$ is arbitrary.
\end{proof}

In view of the previous results, we can apply Theorem~\ref{unique.entropy} to obtain that
$U$ is the unique solution in the sense of Definition \ref{def.entropy} of  \eqref{eq.main} if $q=1+\frac1\gamma$, and   of \eqref{conservation.law} if $q<1+\frac1\gamma$, respectively.
   Since $U$ is uniquely determined,   the whole family $u_\lambda$ converges to it, and not only a subsequence. Moreover,  the uniqueness 
   shows that $U$ is self-similar, $U(t,x)=t^{-\alpha}U_\gamma(xt^{-\beta})$ with $U_\gamma:=U(1,\cdot)$.
 Finally,
  taking $t=1$ and $
 \lambda=t^\beta$, so that $u_\lambda(1,\xi)=t^\alpha u(t,t^\beta\xi)$, we obtain \eqref{loc.conv}.

When $q=1+\frac1\gamma$ we can also prove uniform convergence on compact sets. Here $\eps_\lambda=1$ and, by Lemma~\ref{uniform.gamma.subcritical}(a), $|u_\lambda|\le M$ in $[\frac12,1]\times\rr$, with $M$ independent of $\lambda$. Using
\[
    u_\lambda(1)=G_{1/2}\ast u_\lambda(\tfrac12)-\int_{1/2}^1\partial_xG_{1-s}\ast F\big(u_\lambda(s)\big)\,ds
\]
together with $\|\partial_xG_\tau\|_{L^1(\rr)}\le C\tau^{-1/2}$ and $\|\partial_x^2G_\tau\|_{L^1(\rr)}\le C\tau^{-1}$, we get for $0<|h|\le1$ and $\xi\in\rr$ the following estimate:
\[
    |u_\lambda(1,\xi+h)-u_\lambda(1,\xi)|\le CM|h|+CM^q\int_0^{1/2}\min\Big\{\frac1{\tau^{1/2}},\frac{|h|}\tau\Big\}\,d\tau\le C_M\,|h|\Big(1+\log\frac1{|h|}\Big),
\]
with $C_M$ independent of $\lambda$ and $\xi$. Hence the functions $u_{t^\beta}(1)$, $t>1$, are uniformly bounded and equicontinuous, and, by the Arzel\`a--Ascoli theorem and \eqref{loc.conv}, $U_\gamma$ is continuous and $u_{t^\beta}(1)\to U_\gamma$ uniformly on compact sets as $t\to\infty$.

\subsection{Convergence in the norm of the datum}\label{subsect.upgrade}
We  write $v(t,\xi):=t^{\gamma\theta}u(t,t^\theta\xi)$. The following proposition completes the proof of Theorem~\ref{first.term.gamma<1},
since  $v(t)-U_\gamma=t^{\gamma\theta}\big(u(t)-U(t)\big)(t^\theta\,\cdot)$.

\begin{proposition}\label{ff.global}
Let $\gamma\in(0,1)$, $q>1$, $k_\pm\in\rr$, and let $u_0$ satisfy \eqref{initial.gamma}. Then
 $v(t)\to U_\gamma$ as $t\to\infty$ in $L^{1/\gamma,\infty}(\rr)$ and in $L^p(\rr)$ for every $p\in(\frac1\gamma,\infty)$, and also in $L^\infty(\rr)$ if $q\ge1+\frac1\gamma$.
\end{proposition}

\begin{proof}
Using the far-field estimate \eqref{thm.env} we obtain
\begin{equation}
  \label{eq:est.1}
  |v(t,\xi)|\le\mathcal N\,t^{\gamma\theta}\big(1+t^\theta+t^\theta|\xi|\big)^{-\gamma}\le\mathcal N(1+|\xi|)^{-\gamma}.
\end{equation}
By Theorem~\ref{farfield}, given $\eps\in(0,1]$ there is $L_\eps\ge1$ such that $\big|\,|x|^\gamma u(t,x)-k_\pm\big|\le\eps$ for $|x|\ge L_\eps(1+t^\theta)$. This gives
\begin{equation}
  \label{eq:est.2}
    \big|\,|\xi|^\gamma v(t,\xi)-k_{\sgn\xi}\big|\le\eps\qquad\text{for }t\ge1,\ |\xi|\ge2L_\eps.
\end{equation}
Since $v(t)\to U_\gamma$ in $L^1_{\textup{loc}}(\rr)$ by \eqref{loc.conv}, both bounds  \eqref{eq:est.1} and \eqref{eq:est.2} hold also for $U_\gamma$, a.e., and therefore
\[
    |v(t,\xi)-U_\gamma(\xi)|\le2\eps|\xi|^{-\gamma}\qquad\text{for a.e.\ }|\xi|\ge2L_\eps\text{ and all }t\ge1.
\]

Let $X$ be one of the three spaces.  For $L\ge2L_\eps$ we write $v(t)-U_\gamma=h(t)+\tilde h(t)$, with $h(t):=(v(t)-U_\gamma)\mathds 1_{\{|\xi|<L\}}$. By \eqref{loc.conv} and \eqref{eq:est.1}, $h(t)\to0$ in $L^1(\rr)$  and $|h(t)|\le2\mathcal N$; moreover, if $q\geq 1+ 1/\gamma$, then $h(t)\to 0$ uniformly, by Subsections~\ref{subsect.supercritical} and \ref{subsect.critical}, since $h(t)$ vanishes outside $(-L,L)$. Hence
 $h(t)\to0$ in $X$. For $t\ge1$, $|\tilde h(t)|\le2\eps|\xi|^{-\gamma}\mathds 1_{\{|\xi|\ge L\}}$ a.e., so that \[\|\tilde h(t)\|_X\le2\eps\,\big\||\xi|^{-\gamma}\mathds 1_{\{|\xi|\ge1\}}\big\|_X.\] Hence
\[
    \limsup_{t\to\infty}\|v(t)-U_\gamma\|_X\le C_X\eps,
\]
with $C_X$ depending only on $X$ and $\gamma$, and $\eps\in(0,1]$ is arbitrary.
 \end{proof}

\section{The case \texorpdfstring{$\gamma=1$}{gamma=1}: proof of Theorem~\ref{first.term.gamma=1}}\label{sect-eq1}

Throughout this section $u_0$ satisfies \eqref{initial.gamma} with $\gamma=1$, $g(\lambda):=\lambda/\log\lambda$ for $\lambda>1$,  and
\[\kappa:=\esssup_{x\in\rr}(1+|x|)|u_0(x)|.\]
As in Section~\ref{sect-lt1}, we first prove the convergence in the inner region,
\begin{align}
\label{loc.conv1}
&\lim_{t\to\infty}\int_{|\xi|\le R}\Big|\frac{2t^{1/2}}{\log t}\,u(t,t^{1/2}\xi)-(k_++k_-)G_1(\xi)\Big|\,d\xi=0&&\text{if }q>2,\\
\label{loc.conv2}
&\lim_{\lambda\to\infty}\int_{|\xi|\le R}\big|g(\lambda)\,u\big(\lambda g^{q-1}(\lambda),\lambda\xi\big)-U_{k_++k_-}(1,\xi)\big|\,d\xi=0&&\text{if }1<q\le2,
\end{align}
for every $R>0$, and then upgrade it to \eqref{thm2.weak}--\eqref{thm2.nwave} by means of the far-field estimate \eqref{env.gamma1} of Theorem~\ref{farfield.gamma1}.

\subsection{The case \texorpdfstring{$q>2$}{q>2}}\label{subsect.heat1}
As in Section~\ref{sect-lt1}, we first record the behaviour of the heat equation with data \eqref{initial.gamma}, now with $\gamma=1$, which produces the logarithm.

\begin{lemma}\label{decay.linear.heat1}
    Let $u_0\in L^\infty(\rr)$ satisfy \eqref{initial.gamma} with $\gamma=1$. Then
    \[
        \lim _{t\rightarrow\infty}\frac{t^{1/2}}{\log t}(G_t\ast u_0)(t^{1/2}\xi)=
        \frac{k_++k_-}2G_1(\xi)
    \]
    uniformly for $\xi$ in compact sets.
\end{lemma}

Estimate \eqref{thm2.decay} holds by Proposition~\ref{cor.decay.gamma}. To prove \eqref{loc.conv1} we use the same arguments as in Subsection \ref{subsect.supercritical}. Let  $p_0\in(1,\frac q2)$ such that   $u_0\in L^{p_0}(\rr)$. Denoting $a:=\frac{q-p_0}{2p_0}$  we have
\begin{align}\label{unif.compact.sets.gamma=1}
\frac{t^{1/2}}{\log t}\|u(t)-G_t\ast u_0\|_{L^\infty(\rr)}&\lesssim t^{-1/2}\frac{\max\{1,\log(1+t),t^{1-a}\}}{\log t}+\frac{t^{1-\frac q{2p_0}}}{\log t}.
\end{align}
 Both tend to $0$ as $t\to\infty$, since $a>1/2$ and $ q>{2p_0}$.  Lemma~\ref{decay.linear.heat1} gives \eqref{loc.conv1}.

\subsection{The convective case \texorpdfstring{$1<q\le2$}{1<q<=2}}
For $\lambda\ge e$ let $u_\lambda(t,x):=g(\lambda)\,u\big(\lambda g^{q-1}(\lambda)t,\lambda x\big)$ and $u_{0,\lambda}(x):=g(\lambda)u_0(\lambda x)$. Then $u_\lambda$ solves
\[
\partial_tu_\lambda+\partial_x(F(u_\lambda))=\eps_\lambda\partial_{xx}u_\lambda,\qquad\eps_\lambda:=\frac{g^{q-1}(\lambda)}\lambda=\frac{\lambda^{q-2}}{(\log\lambda)^{q-1}}\in(0,1],
\]
with initial datum $u_{0,\lambda}$, and $\eps_\lambda\to0$ as $\lambda\to\infty$. For any $\psi\in C^1_c(\rr)$ we have
\begin{equation}\label{data.delta}
\int_\rr u_{0,\lambda}\psi\,dx=\frac1{\log\lambda}\int_\rr u_0(y)\psi(y/\lambda)\,dy\longrightarrow(k_++k_-)\psi(0)\qquad\text{as }\lambda\to\infty .
\end{equation}
This will allow us to identify the initial data in the limit equation.

\begin{lemma}\label{estimates.rescaled.gamma=1}
There is a constant $C$ such that, for all $\lambda\ge e$, $t>0$ and $R>0$,
\begin{align}
\label{l.infty.1}
&\|u_\lambda(t)\|_{L^\infty(\rr)}\le Ct^{-\frac1q}\Big(q+\frac{\log(2+t)}{\log\lambda}\Big)^{\frac1q},\\
\label{u.in.balls}
&\int_{|x|<R}|u_\lambda(t,x)|\,dx\le\frac{2\kappa\log(1+\lambda R)}{\log\lambda}\le2\kappa\big(1+\log(1+R)\big).
\end{align}
\end{lemma}

\begin{proof}
Estimate \eqref{l.infty.1} follows from \eqref{thm2.decay}. By a change of variables and  Corollary~\ref{local.L1}, we get
 \[\int_{|x|<R}|u_\lambda(t,x)|\,dx=\frac{g(\lambda)}\lambda\int_{|y|<\lambda R}|u(\lambda g^{q-1}(\lambda)t,y)|\,dy\le\frac{2\kappa\log(1+\lambda R)}{\log\lambda},\]
  and $\log(1+\lambda R)\le\log\lambda+\log(1+R)$ finishes the proof.
\end{proof}

Using \eqref{l.infty.1}, \eqref{u.in.balls},
the same arguments as in Subsection~\ref{subsect.critical} show that there exists a function 
 $U\in C((0,\infty);L^1_{\textup{loc}}(\rr))$ that satisfies \eqref{ineg.k} and, up to a  subsequence,
\begin{equation}\label{n:everyt-gamma1}
    u_{\lambda_j}(t)\to U(t)\quad\text{in }L^1_{\textup{loc}}(\rr)\qquad\text{for every }t>0 .
\end{equation}
By \eqref{l.infty.1}, $\|U(t)\|_{L^\infty(\rr)}\le Cq^{1/q}t^{-1/q}$, and, by \eqref{u.in.balls} and \eqref{n:everyt-gamma1}, $U\in L^1(\rr)$ and satisfies $\|U(t)\|_{L^1(\rr)}\le2\kappa$ for every $t>0$.

Let $\rho\in C^\infty_c(\rr)$ with $0\le\rho\le1$, $\rho=1$ on $[-1,1]$ and $\operatorname{supp}\rho\subset[-2,2]$, and let $\rho_R:=\rho(\cdot/R)$. Since $u_\lambda\in C([0,\infty);L^1_{\textup{loc}}(\rr))$,
\[
    \int_\rr u_\lambda(t)\rho_R\,dx=\int_\rr u_{0,\lambda}\rho_R\,dx+\int_0^t\!\!\int_\rr F(u_\lambda)\rho_R'\,dx\,ds+\eps_\lambda\int_0^t\!\!\int_\rr u_\lambda\rho_R''\,dx\,ds .
\]
The second term on the right is bounded by 
 $Ct^{1/q}/R$, with $C$ independent of $R$ and $t$, by \eqref{l.infty.1}, \eqref{u.in.balls} and dominated convergence; the last term tends to $0$. By \eqref{data.delta} and \eqref{n:everyt-gamma1},
\begin{equation*}
    \Big|\int_\rr U(t)\rho_R\,dx-(k_++k_-)\Big|\le\frac{Ct^{1/q}}R\qquad\text{for all }t,R>0.
\end{equation*}
Letting $R\to\infty$, $\int_\rr U(t)\,dx=k_++k_-$ for all $t>0$. This also gives  a tail control for the limit $U$:
\begin{equation}
\label{tails.U}
\int_{|x|>2R}|U(t,x)|dx\lesssim \frac{t^{1/q}}R\qquad\text{for all }t,R>0.  
\end{equation}
It is first obtained for nonnegative limit and then for changing sign solutions. 

Arguing as above with $\psi\in C^\infty_c(\rr)$ in place of $\rho_R$, we get $\big|\int_\rr U(t)\psi\,dx-(k_++k_-)\psi(0)\big|\le C_\psi t^{1/q}$. For $\psi\in BC(\rr)$ we approximate $\psi\rho_R$ uniformly by smooth compactly supported functions and use $\|U(t)\|_{L^1(\rr)}\le2\kappa$ and \eqref{tails.U}; this gives $\int_\rr U(t)\psi\,dx\to(k_++k_-)\psi(0)$ as $t\downarrow0$ for every $\psi\in BC(\rr)$. 

The limit point $U$ belongs to $L^\infty((0,\infty);L^1(\rr))\cap L^\infty((\tau,\infty)\times\rr)$ for every $\tau>0$, it satisfies \eqref{ineg.k}, and takes the initial datum $(k_++k_-)\delta_0$ in the above sense; by the uniqueness theorem of Liu and Pierre \cite[Theorem~1.2 and Remark~1.3]{MR735207}, $U=U_{k_++k_-}$, the source-type solution of \cite[Section~2]{MR735207}. Since the limit does not depend on the sequence, $u_\lambda(t)\to U_{k_++k_-}(t)$ in $L^1_{\textup{loc}}(\rr)$ for every $t>0$ as $\lambda\to\infty$; for $t=1$, this is \eqref{loc.conv2}.

\subsection{Convergence in the norm of the datum}\label{subsect.upgrade1}
The following proposition completes the proof of Theorem~\ref{first.term.gamma=1}.

\begin{proposition}\label{upgrade.gamma1}
Under the hypotheses of Theorem~\ref{first.term.gamma=1}, let $M:=k_++k_-$. Then, in $X=L^{1,\infty}(\rr)$, in $X=L^p(\rr)$ for every $p\in(1,\infty)$, and, if $q>2$, also in $X=L^\infty(\rr)$,
\begin{align}
\label{resc.heat}
&\lim_{t\to\infty}\Big\|\frac{2t^{1/2}}{\log t}\,u(t,t^{1/2}\cdot)-MG_1\Big\|_{X}=0&&\text{if }q>2,\\
\label{resc.nwave}
&\lim_{\lambda\to\infty}\big\|g(\lambda)\,u\big(\lambda g^{q-1}(\lambda),\lambda\,\cdot\big)-U_M(1,\cdot)\big\|_{X}=0&&\text{if }1<q\le2.
\end{align}
\end{proposition}

\begin{proof}
We set
\[
\Phi:=MG_1, \quad v(\tau,\xi):=\frac{2\tau^{1/2}}{\log\tau}u(\tau,\tau^{1/2}\xi)\quad \text{if}\quad q>2,\]
 \[\Phi:=U_M(1,\cdot),\quad v(\tau,\xi):=g(\tau)u(t_\tau,\tau\xi), \quad t_\tau:=\tau g^{q-1}(\tau)\quad \text{if} \ 1<q\le2.\]
In both cases, by \eqref{thm2.decay} and \eqref{env.gamma1} there is a constant $C$ such that for all  $\tau\ge 3$
\begin{equation}\label{tail.v}
\|v(\tau)\|_{L^\infty(\rr)}\le C,\qquad |v(\tau,\xi)|\le C\Big(\frac1{|\xi|\log\tau}+\frac1{\xi^2}\Big)\quad\text{for }\xi\ne0 .
\end{equation}
In view of  \eqref{loc.conv1}--\eqref{loc.conv2}, $v(\tau)\to\Phi$ in $L^1_{\textup{loc}}(\rr)$, and   by \eqref{unif.compact.sets.gamma=1} and Lemma~\ref{decay.linear.heat1} uniformly on bounded sets for $q>2$, as $\tau\to\infty$.

Let $X$ be one of the three spaces and $L\geq 1$. We write $v(\tau)-\Phi=h(\tau)+\tilde h(\tau)$, with $h(\tau,\xi):=\big(v(\tau,\xi)-\Phi(\xi)\big)\mathds 1_{\{|\xi|<L\}}$. As in the proof of Proposition~\ref{ff.global}, $h(\tau)\to0$ in $X$. By \eqref{tail.v},
\[
    |\tilde h(\tau,\xi)|\le C\Big(\frac1{|\xi|\log\tau}+\frac1{\xi^2}\Big)\mathds 1_{\{|\xi|\ge L\}}+|\Phi(\xi)|\,\mathds 1_{\{|\xi|\ge L\}} .
\]
Since $|\xi|^{-1}\mathds 1_{\{|\xi|\ge1\}}\in X$ and $\Phi\in L^1(\rr)\cap L^\infty(\rr)$, we get
\[
    \limsup_{\tau\to\infty}\|v(\tau)-\Phi\|_X\le C_X\Big(\big\|\xi^{-2}\mathds 1_{\{|\xi|\ge L\}}\big\|_X+\big\|\Phi\mathds 1_{\{|\xi|\ge L\}}\big\|_X\Big)\xrightarrow[L\to\infty]{}0 ,
\]
with $C_X$ depending only on $X$ and $C$. 
\end{proof}

Let us check that Proposition \ref{upgrade.gamma1} indeed implies 
\eqref{thm2.weak}--\eqref{thm2.nwave} which finishes the proof of Theorem~\ref{first.term.gamma=1}.

Let $M:=k_++k_-$. For  $q>2$, 
the identity
\[
\frac{2t^{1/2}}{\log t}\,\big(u(t)-\mathcal U(t)\big)(t^{1/2}\xi)=\frac{2t^{1/2}}{\log t}\,u(t,t^{1/2}\xi)-MG_1(\xi)
\]
and \eqref{resc.heat} give \eqref{thm2.weak}, and \eqref{thm2.heat} after a change of variables.

Let now $1<q\le2$, let $\lambda$ be large, and put $s:=\lambda g^{q-1}(\lambda)=\lambda^q(\log\lambda)^{1-q}$. By \eqref{nwave}, $g(\lambda)^{-1}U_M(1,x/\lambda)=U_{M_1}(s,x)$ with $M_1:=M\log\lambda$, so that, after a change of variables, \eqref{resc.nwave} says that $\|u(s)-U_{M_1}(s)\|_X=o\big(\|U_{M_1}(s)\|_X\big)$ as $s\to\infty$. Let $M_2:=\frac Mq\log s$. Since $\log s=q\log\lambda-(q-1)\log\log\lambda$, we have $M_1/M_2\to1$, and therefore $\|U_{M_1}(s)-U_{M_2}(s)\|_X=o\big(\|U_{M_2}(s)\|_X\big)$. The triangle inequality gives
\[
\|u(s)-U_{M_2}(s)\|_X\le C_X\big(\|u(s)-U_{M_1}(s)\|_X+\|U_{M_1}(s)-U_{M_2}(s)\|_X\big)=o\big(\|U_{M_2}(s)\|_X\big),
\]
which is \eqref{thm2.weak} and \eqref{thm2.nwave}.

\section{Self-similar entropy solutions}\label{sect-profiles}

In this section we prove Theorem~\ref{self.profiles}. Throughout, $F(u)=|u|^{q-1}u$ with $q>1$, $\gamma\in(0,1)$, $(k_+,k_-)\neq(0,0)$, and $\alpha$, $\beta$ are as in \eqref{alpha.beta}; we shall use the relations
\begin{equation}\label{s:elem}
\alpha=\gamma\beta,\qquad \frac1\alpha=q-1+\frac1\gamma,\qquad \beta-\frac1q=\frac{(q-1)(1-\gamma)}{q(1+(q-1)\gamma)}>0 .
\end{equation}

 We look for $U(t,x)=t^{-\alpha}V(xt^{-\beta})$. Where the profile $V$ is smooth, it solves an ordinary differential equation that can be integrated explicitly, and its solutions form a one-parameter family of \emph{branches}; the behaviour of $U_0$ at $\pm\infty$ selects the two branches that carry the tails. These branches cannot be glued continuously, and $V$ is obtained by joining them through admissible shocks. The only delicate point is admissibility. Since $F$ is convex on $[0,\infty)$ and concave on $(-\infty,0]$, a shock between states of opposite signs is admissible only if it is not too strong, and this is what produces the threshold $\eps(\gamma,q)$.

Since $F$ is odd, we may assume that $k_+>0$, or $k_+=0<k_-$.

\subsection{The profile equation and its branches}\label{subsec:profile.eqs}
Let $V\in C^1(I)$ on an open interval $I$, and let $U(t,x)=t^{-\alpha}V(\xi)$, where $\xi=xt^{-\beta}$.  The function $U$ solves $U_t+(F(U))_x=0$ in $\{xt^{-\beta}\in I\}$ if and only if
\begin{equation}\label{eq:ode}
\big(F'(V)-\beta\xi\big)V'=\alpha V\qquad\text{in }I .
\end{equation}
If $V\neq0$ in $I$, then \eqref{eq:ode} forces $V'\neq0$, so that $V$ is strictly monotone and we may regard $\xi=\xi(V)$ as a function of $V$, with $\frac{d\xi}{dV}=\frac1{V'(\xi)}$; then \eqref{eq:ode} becomes
\begin{equation}\label{eq:lin}
F'(V)-\beta\xi=\alpha V\,\frac{d\xi}{dV} .
\end{equation}
Multiplying it by $|V|^{1/\gamma}/(\alpha V)$ and using \eqref{s:elem}, we get
\[
\frac d{dV}\Big(|V|^{1/\gamma}\xi-q|V|^{1/\alpha}\Big)=0 .
\]
Conversely, if $V\neq 0$ and $|V|^{1/\gamma}\xi-q|V|^{1/\alpha}$ is constant in $I$, differentiating in $\xi$ shows that $V'\neq0$, and the computation can be reversed. Hence \eqref{eq:ode} holds in $I$ if and only if there is $C\in\rr$ such that, in $I$,
\begin{equation}\label{eq:firstint}
|V|^{1/\gamma}\xi-q|V|^{1/\alpha}=C,\qquad\text{that is,}\qquad \xi=\Phi_C(|V(\xi)|),
\end{equation}
where
\[\Phi_C(w):=qw^{q-1}+Cw^{-1/\gamma}, \  \text{for} \ w>0.\]
 The same reduction and first integral appear in \cite{MR4687278}. The smooth pieces of $V$ on which $V\neq0$ are thus, up to sign, the inverse functions of the monotone parts of the functions $\Phi_C$, whose shape is elementary.

\begin{lemma}\label{s:lem:phi}
Let $C\in\rr$. Then
\[
\Phi_C'(w)=w^{-\frac1\gamma-1}\Big(q(q-1)\,w^{1/\alpha}-\frac C\gamma\Big),\qquad w>0 .
\]
Consequently:
\begin{enumerate}
\item[(i)] if $C<0$, then $\Phi_C$ is an increasing bijection of $(0,\infty)$ onto $\rr$, and if $C=0$, onto $(0,\infty)$;
\item[(ii)] if $C>0$, then $\Phi_C$ tends to $+\infty$ at $0$ and at $\infty$, it is decreasing on $(0,w_*(C)]$ and increasing on $[w_*(C),\infty)$, where
\[
w_*(C):=\Big(\frac C{\gamma q(q-1)}\Big)^{\alpha},
\]
and its minimum value $\xi_*(C):=\Phi_C(w_*(C))$ satisfies
\begin{equation}\label{s:sonic}
\beta\xi_*(C)=F'(w_*(C))=q\,w_*(C)^{q-1};
\end{equation}
\item[(iii)] $w_*(C)$ and $\xi_*(C)$ are increasing in $C>0$, and $\Phi_C<\Phi_{C'}$ on $(0,\infty)$ if $C<C'$.
\end{enumerate}
\end{lemma}

\begin{proof}
Only \eqref{s:sonic} needs a comment: $\Phi_C$ solves \eqref{eq:lin} on $(0,\infty)$, that is,
\[
F'(w)-\beta\,\Phi_C(w)=\alpha w\,\Phi_C'(w),\qquad w>0,
\]
and $\Phi_C'(w_*(C))=0$.
\end{proof}

\emph{Which branches carry the tails.} Along a branch on which $V\to0$ as $\xi\to\pm\infty$, the first relation in \eqref{eq:firstint} gives $\xi|V(\xi)|^{1/\gamma}\to C$; hence the tails of $U_0$ force
\[
C=-|k_-|^{1/\gamma}\ \text{ on the left},\qquad C=k_+^{1/\gamma}\ \text{ on the right}.
\]
We therefore use the following branches.

\emph{The left branch.} If $k_-\neq0$, let
\[
\Phi_-:=\Phi_{-|k_-|^{1/\gamma}},\qquad w_-:=\Phi_-^{-1},\qquad V_-:=\sgn(k_-)\,w_-\quad\text{on }\rr .
\]
By Lemma~\ref{s:lem:phi}(i), $w_-$ is smooth and increasing, $V_-$ solves \eqref{eq:ode} in $\rr$, and
\[
|\xi|^\gamma V_-(\xi)\to k_-\quad\text{as }\xi\to-\infty,\qquad w_-(\xi)\sim\Big(\frac\xi q\Big)^{\frac1{q-1}}\quad\text{as }\xi\to+\infty .
\]

\emph{The right branch.} If $k_+>0$, let
\[
\Phi_+:=\Phi_{k_+^{1/\gamma}},\qquad w_*:=w_*\big(k_+^{1/\gamma}\big),\qquad \xi_*:=\xi_*\big(k_+^{1/\gamma}\big),
\]
and for $\xi\ge\xi_*$, let $V_+(\xi)\in(0,w_*]$ be the smallest root of $\Phi_+(w)=\xi$, that is, the inverse of the decreasing part of $\Phi_+$. By Lemma~\ref{s:lem:phi}(ii), $V_+$ is continuous and decreasing on $[\xi_*,\infty)$, it solves \eqref{eq:ode} in $(\xi_*,\infty)$, and $\xi^\gamma V_+(\xi)\to k_+$ as $\xi\to\infty$. Since $\xi_*$ is the minimum of $\Phi_+$, attained only at $w_*$, \eqref{s:sonic} gives
\begin{equation}\label{s:Vplus}
V_+(\xi_*)=w_*,\qquad \beta\,\xi_*=F'(w_*)=F'\big(V_+(\xi_*)\big).
\end{equation}

\emph{The rarefaction branch.} For $C=0$, \eqref{eq:firstint} gives the centred rarefaction wave:
\[
w_0(\xi):=\Big(\frac\xi q\Big)^{\frac1{q-1}},\quad \xi>0,\qquad U(t,x)=t^{-\alpha}w_0(xt^{-\beta})=\Big(\frac x{qt}\Big)^{\frac1{q-1}} .
\]

Two branches with different constants never take the same nonzero value at the same point, since by \eqref{eq:firstint} such a value determines $C$. The tails can therefore be connected only through shocks, or through the value $0$.

\subsection{Admissible jumps}\label{subsect.adm}
Let $V$ be bounded, smooth and a solution of \eqref{eq:ode} outside a finite set, and let $\xi_0$ be a point of this set at which $V$ jumps from $v_l=V(\xi_0-)$ to $v_r=V(\xi_0+)\neq v_l$.

The Rankine--Hugoniot condition in similarity variables implies that, if $V$ has a jump at $\xi_0$ with traces $v_l$ and $v_r$, then the shock curve $x/t^\beta=\xi_0 $ satisfies
\begin{equation}\label{RH}
\beta\xi_0=\frac{F(v_r)-F(v_l)}{v_r-v_l}.
\end{equation}
	 Kru\v{z}kov's entropy inequalities are equivalent to the fact (see \cite[Th.~5.1, p.~24]{MR1927887}) that, along the shock,
	\[
	\beta \xi_0 (|v_r-k|-|v_l-k|)\geq \sgn(v_r-k)(F(v_r)-F(k))-\sgn(v_l-k)(F(v_l)-F(k)), \ \forall\ k\in \rr,
	\]
	or, equivalently (see \cite[Lemma 6.1]{MR1304494}), to the Ole\u{\i}nik E-condition along the shock, that is, for all $v$ between $v_l$ and $v_r$,
\begin{equation}\label{eq:oleinik}
\frac{F(v_r)-F(v)}{v_r-v} \leq  \beta \xi_0 \leq \frac{F(v)-F(v_l)}{v-v_l}.
\end{equation}

 Both conditions have a simple geometric meaning. By \eqref{RH}, the speed $s:=\beta\xi_0$ is the slope of the chord $L(v):=F(v_l)+s(v-v_l)$ of the graph of $F$ between $v_l$ and $v_r$, and Ole\u{\i}nik's condition says that the graph of $F$ lies above this chord if $v_l<v_r$ and below it if $v_l>v_r$:
\begin{equation*}
F\ge L\ \text{ between }v_l\text{ and }v_r\ \text{ if }v_l<v_r,\qquad F\le L\ \text{ between }v_r\text{ and }v_l\ \text{ if }v_l>v_r .
\end{equation*}
For our flux, which is convex for $u>0$ and concave for $u<0$, this condition is easily checked.

\begin{lemma}\label{s:lem:adm}
Consider a jump from $v_l$ to $v_r$, and let $s:=\frac{F(v_r)-F(v_l)}{v_r-v_l}$ be the slope of the chord, so that \eqref{RH} reads $\beta\xi_0=s$.
\begin{enumerate}
\item[(i)] If $v_l>v_r\ge0$, or $v_l<v_r\le0$, and \eqref{RH} holds, the jump satisfies \eqref{eq:oleinik}.
\item[(ii)] Let $v_l=-m<0<v_r=p$. Then
\[
s\ge F'(p)\iff p\le r(q)\,m,\qquad s=F'(p)\iff p=r(q)\,m,
\]
where $r(q)$ is the root of \eqref{r.q}. If \eqref{RH} holds, the jump satisfies \eqref{eq:oleinik} if and only if $s\ge F'(p)$.
\end{enumerate}
\end{lemma}

The first property holds by using the convexity of $F$ on $[0,\infty)$ and its concavity on $(-\infty,0]$, while for the second part elementary computations give the desired property.

By \eqref{r.q}, using $\frac1\beta=1+\gamma(q-1)$ and $\gamma<1$, we also get
\begin{equation}\label{s:rident}
r^{1-q}-\frac1\beta=(q-1)(1+r-\gamma)>0,\qquad r=r(q).
\end{equation}

\subsection{Construction of the profiles}\label{subsect.construction}
The choice of the branches below is suggested by the behaviour of $U_0$ at $\pm\infty$; the profiles are then verified directly, and their uniqueness follows from Theorem~\ref{unique.entropy}. To locate a shock between a branch $V_L$ on its left and a branch $V_R$ on its right, with $V_L\neq V_R$ on an interval $J$, we look at
\begin{equation}\label{s:H}
H(\xi):=\beta\xi-S(\xi),\qquad S(\xi):=\frac{F(V_L(\xi))-F(V_R(\xi))}{V_L(\xi)-V_R(\xi)} ,
\end{equation}
the difference between the speed of the similarity line and the slope of the chord; the zeros of $H$ are exactly the points where the two branches can be joined by a jump satisfying \eqref{RH}. The following observation is used three times, in Cases~I, III and IV below, with $C<0$, $C>0$ and $C=0$, respectively: an increasing branch lying above $V_+$ can always be joined to $V_+$ by an admissible shock.

\begin{lemma}\label{s:lem:join}
Let $k_+>0$ and $C<k_+^{1/\gamma}$, and, for $\xi\ge\xi_*$, let $W(\xi)$ be the unique (when $C\leq 0$) or  the largest root (when $C>0$) of $\Phi_C(w)=\xi$. Then $W>V_+$ on $[\xi_*,\infty)$, and there is $\xi_s>\xi_*$ at which the jump from $W(\xi_s)$ to $V_+(\xi_s)$ satisfies \eqref{RH} and \eqref{eq:oleinik}.
\end{lemma}

\begin{proof}
If $C\le0$, then, by Lemma~\ref{s:lem:phi}(i), $W=\Phi_C^{-1}$ is defined on $[\xi_*,\infty)$, since $\xi_*>0$. If $0<C<k_+^{1/\gamma}$, then $\xi_*(C)<\xi_*$, so that, by Lemma~\ref{s:lem:phi}(ii), $W$ defined on $[\xi_*,\infty)$ is the inverse of the restriction of $\Phi_C$ to $[w_*(C),\infty)$. In both cases $W$ is continuous and increasing on $[\xi_*,\infty)$, it solves \eqref{eq:ode} in $(\xi_*,\infty)$, and, since $W(\xi)\to\infty$ and $\xi=qW^{q-1}+CW^{-1/\gamma}$,
\[
W(\xi)\sim\Big(\frac\xi q\Big)^{\frac1{q-1}}\qquad\text{as }\xi\to\infty .
\] Since $C\neq k_+^{1/\gamma}$, the branches $W$ and $V_+$ never take the same value, by \eqref{eq:firstint}, and since $W(\xi)\to\infty$ while $V_+(\xi)\to0$ as $\xi\to\infty$, we get $W>V_+$ on $[\xi_*,\infty)$.

Let $H$ be given by \eqref{s:H} with $V_L=W$, $V_R=V_+$ and $J=[\xi_*,\infty)$. At $\xi_*$ we have $W(\xi_*)>V_+(\xi_*)=w_*$ and $\beta\xi_*=F'(w_*)$, by \eqref{s:Vplus}, so that the strict convexity of $F$ on $[0,\infty)$ gives
\[
S(\xi_*)=\frac{F(W(\xi_*))-F(w_*)}{W(\xi_*)-w_*}>F'(w_*)=\beta\xi_*,\qquad\text{that is,}\qquad H(\xi_*)<0 .
\]
As $\xi\to\infty$, $V_+(\xi)\to0$ while $W(\xi)\sim(\xi/q)^{1/(q-1)}$, so that, by \eqref{s:elem},
\[
S(\xi)\sim\frac\xi q,\qquad H(\xi)\sim\Big(\beta-\frac1q\Big)\xi\to\infty .
\]
Hence $H(\xi_s)=0$ for some $\xi_s>\xi_*$, and the jump at $\xi_s$ satisfies \eqref{eq:oleinik} by Lemma~\ref{s:lem:adm}(i), since $W(\xi_s)>V_+(\xi_s)>0$.
\end{proof}

\begin{figure}[t]
\centering
\includegraphics[width=\linewidth]{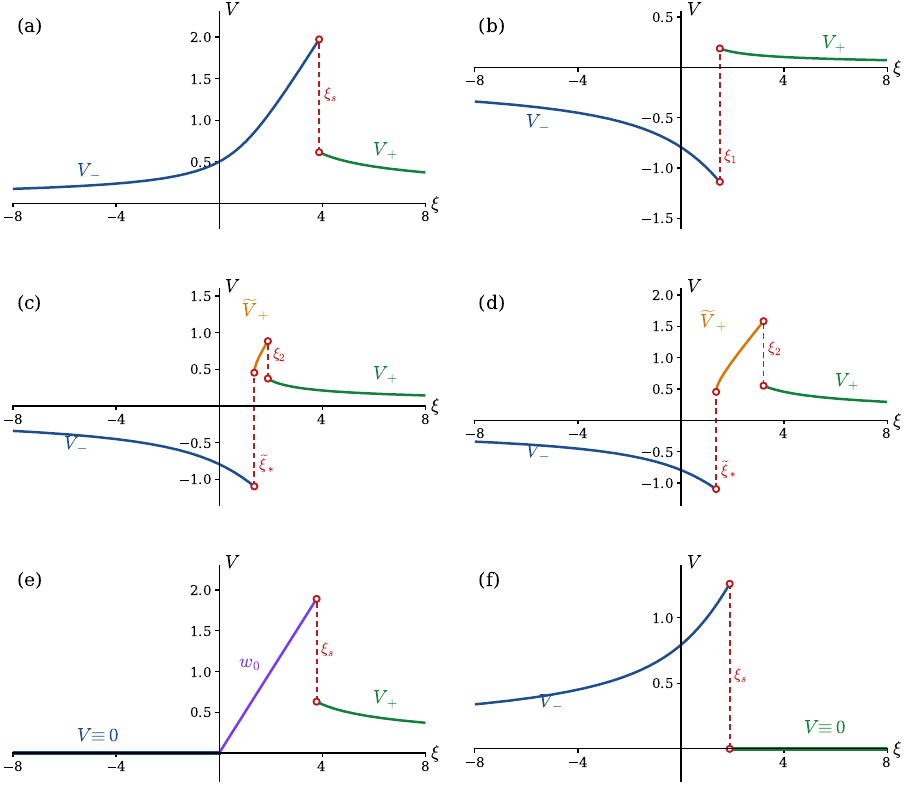}
\caption{Profiles for $q=2$ and $\gamma=\frac12$, for which $\eps(\frac12,2)\approx0.306$. (a)~Case~I, $k_-=0.5$, $k_+=1$. (b)~Case~II, $k_-=-1$, $k_+=0.2$. (c),~(d)~Case~III, $k_-=-1$, and $k_+=0.4$, $k_+=0.8$; the first shock, at $\widetilde\xi_*$, depends only on $k_-$. (e)~Case~IV, $k_-=0$, $k_+=1$; the profile follows $w_0(\xi)=\frac\xi2$ between the origin and the shock. (f)~Case~V, $k_-=1$, $k_+=0$. The notation is that of Subsection~\ref{subsect.construction}.}
\label{fig:profiles}
\end{figure}

\textbf{Case I: $k_+>0$, $k_->0$.} Both tails are positive, and $V_-=w_-$ is the increasing branch with $C=-k_-^{1/\gamma}<k_+^{1/\gamma}$. Applying Lemma~\ref{s:lem:join} with $W=w_-$
 gives $\xi_s>\xi_*$ such that
\[
V(\xi)=\begin{cases}w_-(\xi),&\xi<\xi_s,\\ V_+(\xi),&\xi>\xi_s,\end{cases}
\]
has an admissible jump at $\xi_s$.
See Figure~\ref{fig:profiles}(a).

\textbf{Case II: $k_-<0<k_+\le\eps(\gamma,q)|k_-|$.} Now $V_-=-w_-$ is negative, and a shock from $V_-$ to $V_+$ changes sign; by Lemma~\ref{s:lem:adm}(ii), it is admissible if and only if, at the shock,
\[
\rho(\xi):=\frac{V_+(\xi)}{w_-(\xi)}\le r(q),\qquad \xi\ge\xi_* .
\]
Since $V_+$ decreases and $w_-$ increases, $\rho$ is decreasing on $J:=[\xi_*,\infty)$ and tends to $0$ at infinity, so that we only have to check the condition at $\xi_*$, where, by \eqref{s:Vplus},
\[
V_+(\xi_*)=w_*,\qquad \beta\xi_*=F'(w_*),\qquad \rho(\xi_*)=\frac{w_*}{w_-(\xi_*)} .
\]

It remains to prove that $\rho(\xi_*)\le r(q)$ and that $H(\xi_1)=0$ for some $\xi_1\ge\xi_*$, where $H$ is given by \eqref{s:H} with $V_L=V_-$ and $V_R=V_+$. Then we set $V:=V_-$ on $(-\infty,\xi_1)$ and $V:=V_+$ on $(\xi_1,\infty)$; the jump at $\xi_1$ satisfies \eqref{RH}, and, since $\rho(\xi_1)\le\rho(\xi_*)\le r(q)$, it satisfies \eqref{eq:oleinik} by Lemma~\ref{s:lem:adm}(ii). See Figure~\ref{fig:profiles}(b).

Since $F(V_-)=-w_-^q$, the slope of the chord in \eqref{s:H} is
\[
S(\xi)=\frac{F(V_-(\xi))-F(V_+(\xi))}{V_-(\xi)-V_+(\xi)}=\frac{w_-(\xi)^q+V_+(\xi)^q}{w_-(\xi)+V_+(\xi)},
\]
and $\rho(\xi)\to0$, $w_-(\xi)^{q-1}\sim\frac\xi q$ as $\xi\to\infty$.
As $\xi\to\infty$,
\[
S(\xi)=w_-^{q-1}\,\frac{1+\rho^q}{1+\rho}\sim\frac\xi q,\qquad H(\xi)\sim\Big(\beta-\frac1q\Big)\xi\to\infty .
\] At $\xi_*$ we have $H(\xi_*)=F'(w_*)-S(\xi_*)$, and Lemma~\ref{s:lem:adm}(ii), with $m=w_-(\xi_*)$ and $p=w_*$, shows that $H(\xi_*)\le0$ if and only if $\rho(\xi_*)\le r(q)$. Let $r:=r(q)$. Since $\Phi_-$ is increasing, $\xi_*=\frac q\beta w_*^{q-1}$ by \eqref{s:sonic}, and $w_*^{1/\alpha}=\frac{k_+^{1/\gamma}}{\gamma q(q-1)}$, we get, using \eqref{s:rident} in the third equivalence,
\[
\begin{aligned}
\rho(\xi_*)\le r
&\iff\Phi_-\Big(\frac{w_*}r\Big)\le\xi_*
\iff q\,w_*^{q-1}\Big(r^{1-q}-\frac1\beta\Big)\le|k_-|^{\frac1\gamma}r^{\frac1\gamma}w_*^{-\frac1\gamma}\\
&\iff k_+^{\frac1\gamma}\le|k_-|^{\frac1\gamma}\,\frac{\gamma\,r^{1/\gamma}}{1+r-\gamma}
\iff k_+\le\eps(\gamma,q)\,|k_-| ,
\end{aligned}
\]
and equality holds in all these relations simultaneously. Hence, since $k_+\le\eps(\gamma,q)|k_-|$, we have $\rho(\xi_*)\le r(q)$ and $H(\xi_*)\le0$, and, since $H(\xi)\to\infty$, $H(\xi_1)=0$ for some $\xi_1\ge\xi_*$. The resulting profile is
\[
V(\xi)=\begin{cases}-w_-(\xi),&\xi<\xi_1,\\ V_+(\xi),&\xi>\xi_1;\end{cases}
\]
it decreases from $0$ to $-w_-(\xi_1)$ on $(-\infty,\xi_1)$, jumps up across $0$ to $V_+(\xi_1)$, and decreases to $0$ on $(\xi_1,\infty)$, so that it has exactly one discontinuity.

\textbf{Case III: $k_-<0<k_+$ and $k_+>\eps(\gamma,q)|k_-|$.} Now the shock from $V_-$ to $V_+$ at the  point $\xi_*$ would be too strong, $\rho(\xi_*)>r(q)$, and we first join $V_-$ to an intermediate branch. Let
\[
\widetilde C:=\big(\eps(\gamma,q)|k_-|\big)^{1/\gamma}\in\big(0,k_+^{1/\gamma}\big),\qquad \widetilde w_*:=w_*(\widetilde C),\qquad \widetilde\xi_*:=\xi_*(\widetilde C),
\]
so that $\widetilde\xi_*<\xi_*$ by Lemma~\ref{s:lem:phi}(iii), and, for $\xi\ge\widetilde\xi_*$, let $\widetilde V_+(\xi)$ be the largest root of $\Phi_{\widetilde C}(w)=\xi$, the inverse of the increasing part of $\Phi_{\widetilde C}$; $\widetilde V_+$ is increasing, it solves \eqref{eq:ode} in $(\widetilde\xi_*,\infty)$, and $\widetilde V_+(\widetilde\xi_*)=\widetilde w_*$.

\emph{First shock.} Let $r:=r(q)$. The computation of Case~II, with $k_+$ replaced by $\eps(\gamma,q)|k_-|$, that is, in the case of equality, gives
\[
\Phi_-\Big(\frac{\widetilde w_*}r\Big)=\frac q\beta\,\widetilde w_*^{\,q-1}=\widetilde\xi_*,\qquad\text{that is,}\qquad w_-(\widetilde\xi_*)=\frac{\widetilde w_*}r .
\]
Hence the jump from $V_-(\widetilde\xi_*)=-\widetilde w_*/r$ to $\widetilde V_+(\widetilde\xi_*)=\widetilde w_*$ has $m=\widetilde w_*/r$ and $p=\widetilde w_*=r\,m$, so that Lemma~\ref{s:lem:adm}(ii) and \eqref{s:sonic} give
\[
\frac{F(\widetilde w_*)-F(-\widetilde w_*/r)}{\widetilde w_*+\widetilde w_*/r}=F'(\widetilde w_*)=\beta\,\widetilde\xi_* .
\]
Thus \eqref{RH} holds, and the jump is admissible by Lemma~\ref{s:lem:adm}(ii).

\emph{Second shock.} Lemma~\ref{s:lem:join} with $C=\widetilde C$, for which $W=\widetilde V_+$ on $[\xi_*,\infty)$, gives $\xi_2>\xi_*$ at which $\widetilde V_+$ and $V_+$ are joined by an admissible shock. The lemma applies since $0<\widetilde C<k_+^{1/\gamma}$, and $W=\widetilde V_+$ on $[\xi_*,\infty)$ because $\widetilde\xi_*<\xi_*$.

The profile is given by
\[
V(\xi)=\begin{cases}-w_-(\xi),&\xi<\widetilde\xi_*,\\ \widetilde V_+(\xi),&\widetilde\xi_*<\xi<\xi_2,\\ V_+(\xi),&\xi>\xi_2.\end{cases}
\]
Here $V$ decreases from $0$ to $-\widetilde w_*/r$ on $(-\infty,\widetilde\xi_*)$, jumps up to $\widetilde w_*$, increases to $\widetilde V_+(\xi_2)$, jumps down to $V_+(\xi_2)$, and decreases to $0$. It has exactly two discontinuities, and $|V|=\widetilde V_+$ is increasing between them.
 See Figure~\ref{fig:profiles}(c),~(d).

\textbf{Case IV: $k_+>0$, $k_-=0$.} The left tail vanishes, and the zero state is followed by the rarefaction branch, which starts from $0$ at $\xi=0$ and is the increasing branch with $C=0<k_+^{1/\gamma}$. Lemma~\ref{s:lem:join} with $W=w_0$ gives $\xi_s>\xi_*$, and we set
\[
V(\xi)=\begin{cases}0,&\xi\le0,\\ w_0(\xi)=\big(\frac\xi q\big)^{\frac1{q-1}},&0<\xi<\xi_s,\\ V_+(\xi),&\xi>\xi_s .\end{cases}
\]
See Figure~\ref{fig:profiles}(e).

\textbf{Case V: $k_+=0<k_-$.} The right tail vanishes, and the positive left branch $V_-=w_-$ is joined to $0$. Let $H$ be given by \eqref{s:H} with $V_L=w_-$, $V_R\equiv0$ and $J=\rr$, so that
\[
H(\xi)=\beta\xi-w_-(\xi)^{q-1}\to-\infty\quad\text{as }\xi\to-\infty,\qquad H(\xi)\sim\Big(\beta-\frac1q\Big)\xi\to\infty\quad\text{as }\xi\to\infty .
\]
Here the slope of the chord between $w_-$ and $0$ is $S=\frac{F(w_-)-F(0)}{w_--0}=w_-^{q-1}$, and $w_-(\xi)\to0$ as $\xi\to-\infty$, while $w_-(\xi)^{q-1}\sim\frac\xi q$ as $\xi\to\infty$. Hence $H(\xi_s)=0$ for some $\xi_s$, and $\xi_s>0$, since $\beta\xi_s=w_-(\xi_s)^{q-1}>0$. We set $V:=w_-$ on $(-\infty,\xi_s)$ and $V:=0$ on $(\xi_s,\infty)$. The jump at $\xi_s$ satisfies \eqref{RH}, and \eqref{eq:oleinik} by Lemma~\ref{s:lem:adm}(i). See Figure~\ref{fig:profiles}(f).

\begin{proof}[Proof of Theorem~\ref{self.profiles}]

Let $V$ be the profile constructed above in the corresponding case, and $U(t,x):=t^{-\alpha}V(xt^{-\beta})$. In each case $V$ is bounded,  and $|\xi|^\gamma V(\xi)\to k_\pm$ as $\xi\to\pm\infty$, by the behaviour of $V_\pm$ (and since $V\equiv0$ near $\pm\infty$ if $k_\pm=0$).

In each case $V$ is smooth and solves \eqref{eq:ode} outside a finite set $\{z_1<\dots<z_n\}$, at each point of which it is either continuous or has a jump satisfying \eqref{RH} and \eqref{eq:oleinik}.
By the equivalence recalled in Subsection~\ref{subsect.adm}, $U$ satisfies \eqref{ineg.k}. Moreover $\|U(t)\|_{L^\infty(\rr)}=t^{-\alpha}\|V\|_{L^\infty(\rr)}$, so that $U\in L^\infty((\tau,\infty)\times\rr)$ for every $\tau>0$ and $\sup_{t>0}t^\alpha\|U(t)\|_{L^\infty(\rr)}<\infty$.

 Since $U_0$ is homogeneous of degree $-\gamma$, we have, for $R>0$ and $\lambda:=Rt^{-\beta}$,
\[
\int_{|x|<R}|U(t,x)-U_0(x)|\,dx=t^{\beta-\alpha}\int_{|y|<Rt^{-\beta}}|V-U_0|\,dy
=R^{1-\gamma}\lambda^{\gamma-1}\int_{|y|<\lambda}|V-U_0|\,dy .
\]
Since $V-U_0\in L^1_{\textup{loc}}(\rr)$ and $|y|^\gamma|V(y)-U_0(y)|\to0$ as $|y|\to\infty$, the last integral is $o(\lambda^{1-\gamma})$ as $\lambda\to\infty$, and \eqref{initial.entropy} follows by letting $t\downarrow0$.

Thus $U$ is an entropy solution of \eqref{conservation.law} with initial datum $U_0$ and
\[
\sup_{t>0}t^\alpha\|U(t)\|_{L^\infty(\rr)}<\infty,
\]
and it is unique in this class by Theorem~\ref{unique.entropy}, since $\alpha<\frac1q$.
\end{proof}

\appendix
\section{Convergence at every fixed time}\label{sect-fixedtime}

The compactness arguments of Sections~\ref{sect-lt1} and \ref{sect-eq1} give convergence of the rescaled solutions in $L^1_{\textup{loc}}(Q)$, hence in $L^1_{\textup{loc}}(\rr)$ only for almost every time, whereas \eqref{loc.conv} for $q\le1+\frac1\gamma$ and \eqref{loc.conv2} are statements at a fixed time. This passage from almost every time to every time appears, in different forms, in many works on the large-time behaviour of convection--diffusion equations; see, for instance, \cite{EVZArma,EVZIndiana,Carpio1996}, where it relies on bounds in $W^{1,1}$ obtained from one-sided estimates for nonnegative solutions, or on the $L^1$ contraction for integrable data, and is sometimes only sketched; it seems to be part of the folklore. Since none of these tools is available for sign-changing, non-integrable data, we include a complete argument, based on an $L^1$-comparison between a solution and its time translates.

In this appendix $F\in C^1(\rr)$ is an arbitrary flux with $F(0)=0$ (in Sections~\ref{sect-lt1} and \ref{sect-eq1} the results are applied with $F(u)=|u|^{q-1}u$); for $M>0$ we set $L_M:=\max_{[-M,M]}|F'|$, and we let $\omega(x):=e^{-\langle x\rangle}$, $\langle x\rangle:=(1+x^2)^{1/2}$, $L^1_\omega:=L^1(\rr,\omega\,dx)$ and $\|g\|_\omega:=\int_\rr|g|\,\omega\,dx$. Since $\langle x\rangle'\in[-1,1]$ and $\langle x\rangle''\in(0,1]$, we have $\omega\in L^1(\rr)$, $|\omega'|\le \omega$ and $|\omega''|\le2\omega$.

\begin{lemma}\label{n:comparison}
Let $0<\tau<T$, $M>0$, $\eps\in(0,1]$, and let $v_1,v_2$ be classical solutions of
\begin{equation*}
\partial_tv+\partial_x(F(v))=\eps\,\partial_{xx}v\qquad\text{in }(\tau,T)\times\rr
\end{equation*}
with $|v_1|,|v_2|\le M$. Then, with $C_M:=L_M+2$,
\begin{equation}\label{n:comp}
\|v_1(t)-v_2(t)\|_\omega\le e^{C_M(t-t')}\,\|v_1(t')-v_2(t')\|_\omega\qquad\text{for all }\tau<t'\le t<T .
\end{equation}
\end{lemma}

\begin{proof}
Let $z:=v_1-v_2$, and let $\chi\in C^\infty_c(\rr)$ with $0\le\chi\le1$, $\chi=1$ on $[-1,1]$, and $\chi_n:=\chi(\cdot/n)$. By Lemma~\ref{u:kato} with $\phi=\omega\chi_n$, which uses only the local Lipschitz continuity of $F$, and since $|F(v_1)-F(v_2)|\le L_M|z|$,
\[
\int_\rr|z(t)|\omega\chi_n\,dx\le\int_\rr|z(t')|\omega\chi_n\,dx+\int_{t'}^t\!\!\int_\rr|z|\big(L_M|(\omega\chi_n)'|+\eps|(\omega\chi_n)''|\big)\,dx\,ds .
\]
As $n\to\infty$, $(\omega\chi_n)'\to \omega'$ and $(\omega\chi_n)''\to \omega''$ pointwise, with $|(\omega\chi_n)'|+|(\omega\chi_n)''|\le C\omega$. Since $|z|\le2M$, dominated convergence gives $\|z(t)\|_\omega\le\|z(t')\|_\omega+C_M\int_{t'}^t\|z(s)\|_\omega\,ds$, and \eqref{n:comp} follows from Gronwall's lemma, since $s\mapsto\|z(s)\|_\omega$ is continuous.
\end{proof}

Throughout the rest of this appendix, $0<\tau<T$, $\eps_j\in(0,1]$, and $u_j$ are classical solutions of $\partial_tu_j+\partial_x(F(u_j))=\eps_j\partial_{xx}u_j$ in $(\tau,T)\times\rr$, $j\in\nn$, such that
\begin{enumerate}
\item[\textup{(H1)}] $M:=\sup_j\|u_j\|_{L^\infty((\tau,T)\times\rr)}<\infty$;
\item[\textup{(H2)}] $u_j\to U$ in $L^1_{\textup{loc}}((\tau,T)\times\rr)$ as $j\to\infty$.
\end{enumerate}
Note that $|U|\le M$ a.e., and that, since $\omega\in L^1(\rr)$, \textup{(H1)} and \textup{(H2)} imply
\begin{equation}\label{n:L1time}
\int_{t_1}^{t_2}\|u_j(r)-U(r)\|_\omega\,dr\to0\qquad\text{as }j\to\infty,\quad\text{for all }\tau<t_1<t_2<T .
\end{equation}

\begin{lemma}\label{n:timemod}
Let $\tau<a<b<T$, $\eta:=\min\{1,\frac{a-\tau}3,\frac{T-b}2\}$, and
\[
\Lambda(\sigma):=\frac{e^{C_M\eta}}\eta\,\sup_{0<h\le\sigma}\int_{a-2\eta}^b\|U(r+h)-U(r)\|_\omega\,dr,\qquad0<\sigma\le\eta .
\]
Then $\Lambda$ is nondecreasing, $\Lambda(\sigma)\to0$ as $\sigma\downarrow0$, and
\begin{equation}\label{n:timemod-est}
\limsup_{j\to\infty}\ \sup_{0<h\le\sigma}\ \sup_{a-\eta\le s\le b}\|u_j(s+h)-u_j(s)\|_\omega\le\Lambda(\sigma)\qquad\text{for }0<\sigma\le\eta .
\end{equation}
\end{lemma}

\begin{proof}
Let $0<h\le\eta$. Since the equation is autonomous, $u_j(\cdot+h,\cdot)$ is a classical solution of the same equation in $(\tau,T-h)\times\rr$, bounded by $M$, and $b<T-h$. Lemma~\ref{n:comparison} gives $\|u_j(s+h)-u_j(s)\|_\omega\le e^{C_M(s-r)}\|u_j(r+h)-u_j(r)\|_\omega$ for $\tau<r\le s\le b$. For $a-\eta\le s\le b$ we have $s-\eta\ge a-2\eta>\tau$; averaging over $r\in(s-\eta,s)$ and inserting $U$, we get, for $0<h\le\sigma\le\eta$,
\begin{align*}
\|u_j(s+h)-u_j(s)\|_\omega&\le\frac{e^{C_M\eta}}\eta\int_{s-\eta}^s\|u_j(r+h)-u_j(r)\|_\omega\,dr\\
&\le\Lambda(\sigma)+\frac{2e^{C_M\eta}}\eta\int_{a-2\eta}^{b+\eta}\|u_j(r)-U(r)\|_\omega\,dr .
\end{align*}
By \eqref{n:L1time} this gives \eqref{n:timemod-est}. Finally, extending $U\omega$ by zero outside $(a-2\eta,b+\eta)\times\rr$ we obtain a function in $L^1(\rr^2)$, and the continuity of translations in $L^1(\rr^2)$ shows that $\Lambda(\sigma)\to0$ as $\sigma\downarrow0$.
\end{proof}

\begin{proposition}\label{n:fixedtime}
Assume \textup{(H1)} and \textup{(H2)}. Then $U$ has a representative $\widetilde U\in C((\tau,T);\allowbreak L^1_{\textup{loc}}(\rr))$, and $u_j(t)\to\widetilde U(t)$ in $L^1_{\textup{loc}}(\rr)$ for every $t\in(\tau,T)$.
\end{proposition}

\begin{proof}
Since $\omega\ge e^{-\langle R\rangle}$ on $[-R,R]$, it suffices to prove the statement in $L^1_\omega$ for $t$ in a compact interval $[a,b]\subset(\tau,T)$. Let $\eta$ and $\Lambda$ be as in Lemma~\ref{n:timemod}. By \eqref{n:L1time}, for a subsequence, not relabelled, there is a set $E\subset[a-2\eta,b+\eta]$ of full measure such that
\begin{equation}\label{n:goodtimes}
\|u_j(s)-U(s)\|_\omega\to0\qquad\text{for every }s\in E .
\end{equation}

 Let $t\in[a,b]$ and $\varrho>0$. Choose $\sigma\in(0,\eta]$ with $\Lambda(\sigma)\le\varrho$ and $s\in E\cap[t-\sigma,t)$. Since $\|u_j(t)-U(s)\|_\omega\le\|u_j(s+(t-s))-u_j(s)\|_\omega+\|u_j(s)-U(s)\|_\omega$, \eqref{n:timemod-est} and \eqref{n:goodtimes} give $\limsup_j\|u_j(t)-U(s)\|_\omega\le\varrho$. Hence $(u_j(t))_j$ is a Cauchy sequence in $L^1_\omega$; let $\widetilde U(t)$ be its limit. By \eqref{n:goodtimes}, $\widetilde U(t)=U(t)$ for $t\in E\cap[a,b]$, so $\widetilde U=U$ a.e.\ in $(a,b)\times\rr$.

 By the above convergence for each fixed $t$ and \eqref{n:timemod-est}, 
 \[
 \|\widetilde U(t)-\widetilde U(s)\|_\omega=\lim_j\|u_j(t)-u_j(s)\|_\omega\le\Lambda(t-s)\]
  for $a\le s<t\le b$ with $t-s\le\eta$; hence $\widetilde U\in C([a,b];L^1_\omega)$.

Every subsequence of $(u_j)$ satisfies \textup{(H1)} and \textup{(H2)}, and hence, as we have seen, has a further subsequence converging in $L^1_\omega$, for every $t\in[a,b]$, to a continuous function from $[a,b]$ into $L^1_\omega$ which equals $U(t)$ for a.e.\ $t$. Two such functions agree everywhere, so the limit $\widetilde U(t)$ is the same for all subsequences, and the whole sequence converges for every $t\in[a,b]$. The representatives obtained on different intervals $[a,b]$ agree on their intersections, which concludes the proof.
\end{proof}

\section*{Acknowledgements}

N.~De Nitti is a member of the Gruppo Nazionale per l'Analisi Matematica, la Probabilit\`a e le loro Applicazioni (GNAMPA) of the Istituto Nazionale di Alta Matematica (INdAM) and has received support from the INdAM--GNAMPA Project 2026 \textit{Modelli Non-locali in Fluidodinamica, Traffico ed Elasticit\`a} (CUP:~E53C25002010001).

L.~I.~Ignat  was partially supported by a grant of the Ministry of Research, Innovation, and Digitization, CCCDI -
UEFISCDI, project number ROSUA-2024-0001, within PNCDI IV.

{\sloppy F.~Quir\'os was supported by grants PID2023-146931NB-I00, RED2024-153842-T and CEX2023-001347-S, all of them funded by MICIU/AEI/10.13039/501100011033.\par}

This article is based upon work from COST Action 24122 mSPACE, supported by COST (European Cooperation in Science and Technology), www.cost.eu.

\section*{AI Disclosure}

The authors acknowledge the use of  Anthropic Claude as an AI tool for proofreading and for bibliographic verification. The authors independently checked the references against the original sources. All mathematical arguments and proofs in the final manuscript were checked and written by the authors, who retain full responsibility for the whole content.

\bibliographystyle{abbrv}
\bibliography{biblio}

\begin{thebibliography}{10}

\bibitem{AlBe22}
D.~Albritton and R.~Beekie.
\newblock Long-time behavior of scalar conservation laws with critical dissipation.
\newblock {\em Ann. Inst. H. Poincar\'e C Anal. Non Lin\'eaire}, 39(1):225--243, 2022.

\bibitem{AlRo84}
N.~D. Alikakos and R.~Rostamian.
\newblock On the uniformization of the solutions of the porous medium equation in {${\bf R}^N$}.
\newblock {\em Israel J. Math.}, 47(4):270--290, 1984.

\bibitem{MR3381494}
B.~Andrews.
\newblock Moduli of continuity, isoperimetric profiles, and multi-point estimates in geometric heat equations.
\newblock In {\em Surveys in differential geometry 2014. {R}egularity and evolution of nonlinear equations}, volume~19 of {\em Surv. Differ. Geom.}, pages 1--47. Int. Press, Somerville, MA, 2015.

\bibitem{MR4294711}
B.~Andrews.
\newblock Multi-point maximum principles and eigenvalue estimates.
\newblock In {\em 2018 {MATRIX} annals}, volume~3 of {\em MATRIX Book Ser.}, pages 185--213. Springer, Cham, 2020.

\bibitem{MR2517770}
B.~Andrews and J.~Clutterbuck.
\newblock Lipschitz bounds for solutions of quasilinear parabolic equations in one space variable.
\newblock {\em J. Differential Equations}, 246(11):4268--4283, 2009.

\bibitem{MR2504416}
B.~Andrews and J.~Clutterbuck.
\newblock Time-interior gradient estimates for quasilinear parabolic equations.
\newblock {\em Indiana Univ. Math. J.}, 58(1):351--380, 2009.

\bibitem{MR1784177}
P.~B\'enilan, J.~Carrillo, and P.~Wittbold.
\newblock Renormalized entropy solutions of scalar conservation laws.
\newblock {\em Ann. Scuola Norm. Sup. Pisa Cl. Sci. (4)}, 29(2):313--327, 2000.

\bibitem{BeSh88}
C.~Bennett and R.~Sharpley.
\newblock {\em Interpolation of operators}, volume 129 of {\em Pure and Applied Mathematics}.
\newblock Academic Press, Boston, MA, 1988.

\bibitem{MR4671449}
L.~Boccardo, L.~Orsina, and M.~M. Porzio.
\newblock Asymptotic behavior of solutions for nonlinear parabolic problems with {M}arcinkiewicz data.
\newblock {\em J. Evol. Equ.}, 23(4):Paper No. 77, 10, 2023.

\bibitem{Carpio1996}
A.~Carpio.
\newblock Large time behaviour in convection-diffusion equations.
\newblock {\em Ann. Scuola Norm. Sup. Pisa Cl. Sci. (4)}, 23(3):551--574, 1996.

\bibitem{CazacuIgnatPazoto}
C.~Cazacu, L.~Ignat, and A.~Pazoto.
\newblock On the asymptotic behavior of a subcritical convection-diffusion equation with nonlocal diffusion.
\newblock {\em Nonlinearity}, 30(8):3126--3150, 2017.

\bibitem{Daf77}
C.~M. Dafermos.
\newblock Characteristics in hyperbolic conservation laws. {A} study of the structure and the asymptotic behaviour of solutions.
\newblock In {\em Nonlinear analysis and mechanics: {H}eriot-{W}att {S}ymposium, {V}ol. {I}}, volume~17 of {\em Res. Notes in Math.}, pages 1--58. Pitman, London, 1977.

\bibitem{Daf2016}
C.~M. Dafermos.
\newblock {\em Hyperbolic conservation laws in continuum physics}, volume 325 of {\em Grundlehren der mathematischen Wissenschaften [Fundamental Principles of Mathematical Sciences]}.
\newblock Springer-Verlag, Berlin, fourth edition, 2016.

\bibitem{2406.10874}
N.~De~Nitti and E.~Pacherie.
\newblock Nested discontinuous asymptotic profiles for the viscous {B}urgers' equation with infinite mass.
\newblock {\em Asymptot. Anal.}, 2026.
\newblock Published online, doi:10.1177/09217134261469833.

\bibitem{Dix2002}
D.~B. Dix.
\newblock Large-time behaviour of solutions of {B}urgers' equation.
\newblock {\em Proc. Roy. Soc. Edinburgh Sect. A}, 132(4):843--878, 2002.

\bibitem{endalJMPA}
J.~Endal, L.~I. Ignat, and F.~Quir\'os.
\newblock Large-time behaviour for anisotropic stable nonlocal diffusion problems with convection.
\newblock {\em J. Math. Pures Appl. (9)}, 179:277--336, 2023.

\bibitem{endalDCDS}
J.~Endal, L.~I. Ignat, and F.~Quir\'os.
\newblock Heat equations with fast convection: source-type solutions and large-time behaviour.
\newblock {\em Discrete Contin. Dyn. Syst. Ser. S}, 17(4):1497--1512, 2024.

\bibitem{endal}
J.~Endal and E.~R. Jakobsen.
\newblock {$L^1$} contraction for bounded (nonintegrable) solutions of degenerate parabolic equations.
\newblock {\em SIAM J. Math. Anal.}, 46(6):3957--3982, 2014.

\bibitem{EVZArma}
M.~Escobedo, J.~L. V{\'a}zquez, and E.~Zuazua.
\newblock Asymptotic behaviour and source-type solutions for a diffusion-convection equation.
\newblock {\em Arch. Rational Mech. Anal.}, 124(1):43--65, 1993.

\bibitem{EVZIndiana}
M.~Escobedo, J.~L. V\'azquez, and E.~Zuazua.
\newblock A diffusion-convection equation in several space dimensions.
\newblock {\em Indiana Univ. Math. J.}, 42(4):1413--1440, 1993.

\bibitem{EZ}
M.~Escobedo and E.~Zuazua.
\newblock Large time behavior for convection-diffusion equations in {${\bf R}^N$}.
\newblock {\em J. Funct. Anal.}, 100(1):119--161, 1991.

\bibitem{EsZu97}
M.~Escobedo and E.~Zuazua.
\newblock Long-time behavior for a convection-diffusion equation in higher dimensions.
\newblock {\em SIAM J. Math. Anal.}, 28(3):570--594, 1997.

\bibitem{MR4687278}
T.-E. Ghoul, N.~Masmoudi, and E.~Pacherie.
\newblock Nonlinear enhanced dissipation in viscous {B}urgers type equations {II}.
\newblock {\em Comm. Partial Differential Equations}, 48(10-12):1286--1322, 2023.

\bibitem{MR4899801}
T.-E. Ghoul, N.~Masmoudi, and E.~Pacherie.
\newblock Nonlinear enhanced dissipation in viscous {B}urgers-type equations.
\newblock {\em Ann. Inst. H. Poincar\'{e} C Anal. Non Lin\'{e}aire}, 42(3):755--805, 2025.

\bibitem{MR1304494}
E.~Godlewski and P.-A. Raviart.
\newblock {\em Hyperbolic systems of conservation laws}, volume 3/4 of {\em Math\'ematiques \& Applications (Paris) [Mathematics and Applications]}.
\newblock Ellipses, Paris, 1991.

\bibitem{MR3243734}
L.~Grafakos.
\newblock {\em Classical {F}ourier analysis}, volume 249 of {\em Graduate Texts in Mathematics}.
\newblock Springer, New York, third edition, 2014.

\bibitem{Herraiz-1999}
L.~Herraiz.
\newblock Asymptotic behaviour of solutions of some semilinear parabolic problems.
\newblock {\em Ann. Inst. H. Poincar\'{e} Anal. Non Lin\'{e}aire}, 16(1):49--105, 1999.

\bibitem{Hof83}
D.~Hoff.
\newblock The sharp form of {O}le\u{\i}nik's entropy condition in several space variables.
\newblock {\em Trans. Amer. Math. Soc.}, 276(2):707--714, 1983.

\bibitem{MR3789847}
L.~I. Ignat and D.~Stan.
\newblock Asymptotic behavior of solutions to fractional diffusion-convection equations.
\newblock {\em J. Lond. Math. Soc. (2)}, 97(2):258--281, 2018.

\bibitem{JaSe23}
T.~Jakubowski and G.~Serafin.
\newblock Fractional {B}urgers equation with singular initial condition.
\newblock {\em Nonlinear Anal.}, 232:Paper No. 113269, 2023.

\bibitem{KaUg87}
S.~Kamin and M.~Ughi.
\newblock On the behaviour as {$t\to\infty$} of the solutions of the {C}auchy problem for certain nonlinear parabolic equations.
\newblock {\em J. Math. Anal. Appl.}, 128(2):456--469, 1987.

\bibitem{KarchMiaoXu}
G.~Karch, C.~Miao, and X.~Xu.
\newblock On convergence of solutions of fractal {B}urgers equation toward rarefaction waves.
\newblock {\em SIAM J. Math. Anal.}, 39(5):1536--1549, 2008.

\bibitem{Kim2004}
Y.~J. Kim.
\newblock An {O}leinik-type estimate for a convection-diffusion equation and convergence to {$N$}-waves.
\newblock {\em J. Differential Equations}, 199(2):269--289, 2004.

\bibitem{Kru67}
S.~N. Kru\v{z}kov.
\newblock Nonlinear parabolic equations with two independent variables.
\newblock {\em Trudy Moskov. Mat. Ob\v{s}\v{c}.}, 16:329--346, 1967.
\newblock (Russian).

\bibitem{MR267257}
S.~N. Kru\v{z}kov.
\newblock First order quasilinear equations with several independent variables.
\newblock {\em Mat. Sb. (N.S.)}, 81(123):228--255, 1970.

\bibitem{LaurencotFast}
P.~Lauren\c{c}ot.
\newblock Long-time behaviour for diffusion equations with fast convection.
\newblock {\em Ann. Mat. Pura Appl. (4)}, 175:233--251, 1998.

\bibitem{MR1927887}
P.~G. LeFloch.
\newblock {\em Hyperbolic systems of conservation laws}.
\newblock Lectures in Mathematics ETH Z\"urich. Birkh\"auser Verlag, Basel, 2002.
\newblock The theory of classical and nonclassical shock waves.

\bibitem{Lie96}
G.~M. Lieberman.
\newblock {\em Second order parabolic differential equations}.
\newblock World Scientific Publishing Co., River Edge, NJ, revised edition, 2005.

\bibitem{MR735207}
T.-P. Liu and M.~Pierre.
\newblock Source-solutions and asymptotic behavior in conservation laws.
\newblock {\em J. Differential Equations}, 51(3):419--441, 1984.

\bibitem{MR2839869}
P.~V. Lysuho and E.~Y. Panov.
\newblock Renormalized entropy solutions to the {Cauchy} problem for a quasilinear first order equation.
\newblock {\em J. Math. Sci. (N.Y.)}, 172(1):1--23, 2011.

\bibitem{Ole63}
O.~A. Ole\u{\i}nik.
\newblock Discontinuous solutions of non-linear differential equations.
\newblock {\em Amer. Math. Soc. Transl. (2)}, 26:95--172, 1963.

\bibitem{zbMATH01799645}
J.~L. V{\'a}zquez and E.~Zuazua.
\newblock Complexity of large time behaviour of evolution equations with bounded data.
\newblock {\em Chinese Ann. Math. Ser. B}, 23(2):293--310, 2002.

\bibitem{watanabe}
H.~Watanabe.
\newblock {$N$}-wave-like properties for entropy solutions to scalar parabolic-hyperbolic conservation laws.
\newblock {\em Nonlinear Anal. Real World Appl.}, 84:Paper No. 104265, 13, 2025.

\bibitem{zuazuabook}
E.~Zuazua.
\newblock Asymptotic behavior of scalar convection-diffusion equations.
\newblock Preprint, Springer, 2025.

\end{thebibliography}

\end{document}